\documentclass[final]{siamart1116}

\usepackage{amsfonts,amssymb,amsmath,amsopn,float}
\usepackage{tikz}
\usepackage{pgfplots}
\pgfplotsset{compat=1.13}
\usepgfplotslibrary{units}
\usepackage{subcaption}

\newcommand{\bbR}{\mathbb{R}}

\newcommand{\dyad}{\mathbf{\otimes}}
\newcommand{\bzero}{\mathbf{0}}

\newcommand{\bxi}{{\boldsymbol{\xi}}}
\newcommand{\btau}{{\boldsymbol{\tau}}}
\newcommand{\bsigma}{{\boldsymbol{\sigma}}}
\newcommand{\betanew}{\boldsymbol{\eta}}

\newcommand{\Ltwo}[1]{%
\ifthenelse{\equal{#1}{}}{L^2}{L^2(#1)}%
}

\newcommand{\Ltwoz}[1]{%
\ifthenelse{\equal{#1}{}}{L^2_0}{L^2_0(#1)}%
}

\newcommand{\Cone}[1]{%
\ifthenelse{\equal{#1}{}}{C^{1}}{C^{1}(#1)}%
}

\newcommand{\Conez}[1]{%
\ifthenelse{\equal{#1}{}}{C^{1}_{0}}{C^{1}_{0}(#1)}%
}

\newcommand{\Ctwo}[1]{%
\ifthenelse{\equal{#1}{}}{C^{2}}{C^2(#1)}%
}

\newcommand{\Ctwoz}[1]{%
\ifthenelse{\equal{#1}{}}{C^{2}_{0}}{C^{2}_{0}(#1)}%
}

\newcommand{\Cholder}[1]{%
\ifthenelse{\equal{#1}{}}{C^{0,\gamma}}{C^{0,\gamma}(#1)}%
}

\newcommand{\Cholderz}[1]{%
\ifthenelse{\equal{#1}{}}{C^{0,\gamma}_{0}}{C^{0,\gamma}_{0}(#1)}%
}

\newcommand{\bolds}[1]{\boldsymbol{#1}}
\newcommand{\ba}{\bolds{a}}
\newcommand{\bb}{\bolds{b}}
\newcommand{\bc}{\bolds{c}}

\newcommand{\be}{\bolds{e}}

\newcommand{\br}{\bolds{r}}

\newcommand{\bu}{\bolds{u}}
\newcommand{\bv}{\bolds{v}}
\newcommand{\bw}{\bolds{w}}
\newcommand{\bx}{\bolds{x}}
\newcommand{\by}{\bolds{y}}

\newcommand{\bubar}{\bar{\bu}}
\newcommand{\butilde}{\tilde{\bu}}

\newcommand{\sref}[2]{\hyperref[#2]{#1 \ref*{#2}}} 

\newcommand{\TheTitle}{Finite element convergence for state-based peridynamic fracture models} 

\newcommand{\TheAuthors}{Prashant K. Jha and Robert Lipton}

\title{{\TheTitle}\thanks{\textbf{Funding: }This material is based upon work supported by the U. S. Army Research Laboratory and the U. S. Army Research Office under contract/grant number W911NF1610456.}}

\author{
  Prashant K. Jha
  \thanks{Department of Mathematics, Louisiana State University, Baton Rouge, LA 
  (\email{prashant.j16o@gmail.com}). Orcid: https://orcid.org/0000-0003-2158-364X}
  \and
  Robert Lipton
  \thanks{Department of Mathematics, Louisiana State University, Baton Rouge, LA 
  (\email{lipton@math.lsu.edu}). Orcid: https://orcid.org/0000-0002-1382-3204}
}

\headers{Finite element convergence for state-based Peridynamic}{\TheAuthors}

\ifpdf
  \DeclareGraphicsExtensions{.eps,.pdf,.png,.jpg}
\else
  \DeclareGraphicsExtensions{.eps}
\fi

\ifpdf
\hypersetup{
  pdftitle={\TheTitle},
  pdfauthor={\TheAuthors}
}
\fi

\begin{document}

\maketitle

\begin{abstract}
We establish the a-priori convergence rate for finite element approximations of a class of nonlocal nonlinear fracture models. We consider state based peridynamic models where the force at a material point is due to both the strain between two points and the change in volume inside the domain of nonlocal interaction. The  pairwise interactions between points are mediated by a bond potential of multi-well type while multi point interactions are associated with volume change mediated by a hydrostatic strain potential. The hydrostatic potential can either be a quadratic function, delivering a linear force-strain relation, or a multi-well type that can be associated with material degradation and cavitation. We first show the well-posedness of the peridynamic formulation and that peridynamic evolutions exist in the Sobolev space $H^2$. We show that the finite element approximations converge to the $H^2$ solutions uniformly as measured in the mean square norm. For linear continuous finite elements the convergence rate is shown to be $C_t \Delta t + C_s h^2/\epsilon^2$, where $\epsilon$ is the size of horizon, $h$ is the mesh size, and $\Delta t$ is the size of time step. The constants $C_t$ and $C_s$ are independent of $\Delta t$ and $h$ and may depend on $\epsilon$ through the norm of the exact solution. We demonstrate the stability of the semi-discrete approximation. The stability of the fully discrete approximation is shown for the linearized peridynamic force. We present numerical simulations with dynamic crack propagation that support the theoretical convergence rate.

\end{abstract}



\section{Introduction}
In this work, we study the state-based peridynamic theory and obtain an a-priori error bound for the finite element approximation. The peridynamic theory is a reformulation of classical continuum mechanics carried out in the work of Silling in \cite{CMPer-Silling,States}. The strain inside the medium is expressed in terms of  displacement differences as opposed to the displacement gradients. Acceleration of a point is now due to the sum of the forces acting on the point from near by points. The new kinematics bypasses the difficulty incurred by juxtaposing displacement gradients and discontinuities as in the case of classical fracture theories. The nonlocal fracture theory has been applied numerically  to model the complex fracture phenomenon in materials, see \cite{WeckAbe}, \cite{SillBob}, \cite{CMPer-Silling4}, \\ \cite{CMPer-Silling5}, \cite{CMPer-Silling7}, \cite{HaBobaru}, \cite{CMPer-Agwai}, \\\cite{BobaruHu}, \cite{CMPer-Ghajari}, \cite{CMPer-Du}, \cite{CMPer-Lipton2}. 
Every point interacts with its neighbors inside a ball of fixed radius called the horizon.  The size of horizon sets the length scale of nonlocal interaction. When the forces between points are linear and nonlocal length scale tends to zero, it is seen that peridynamic models converge to the classic model of linear elasticity, \cite{CMPer-Emmrich,CMPer-Silling4,AksoyluUnlu,CMPer-Mengesha2}.  The work of \cite{Tian-Du} provides an analytic
framework for analyzing FEM for linear bond and state-based peridynamics.
For nonlinear forces associated with double well potentials the peridynamic evolution converges in the small horizon limit to an evolution with a sharp evolving fracture set and the evolution is governed by the classic linear elastic wave equation away from the fracture set, see \cite{CMPer-Lipton,CMPer-Lipton3,CMPer-JhaLiptonFD}. 
A recent review of the state of the art can be found in \cite{Handbook} and \cite{Du 2018a}.

\begin{figure}
    \centering
    \begin{subfigure}{.45\linewidth}
    \centering
        \begin{tikzpicture}[xscale=0.75,yscale=0.75]
		    \draw [<->,thick] (0,5) -- (0,0) -- (3.0,0);
			\draw [-,thick] (0,0) -- (-3.5,0);
			\draw [-,thin] (0,2.15) -- (2.5,2.15);
			\draw [-, thin] (-3.5,4.15) -- (0,4.15);
			\draw [-,thick] (0,0) to [out=0,in=-175] (2.5,2);
			\draw [-,thick] (-3.5,4) to [out=-5,in=180] (0,0);
			
			\draw (1.5,-0.2) -- (1.5, 0.2);
			\node [below] at (1.5,-0.2) {${r}^+$};
			
			\draw (-2.25,-0.2) -- (-2.25, 0.2);
			\node [below] at (-2.0,-0.2) {${r}^-$};

			\node [right] at (3,0) {$r$};
			\node [left] at (0,2.250) {$C^+$};
			\node [left] at (0.9,4.15) {$C^-$};
			\node [above] at (-2.5,2.20) {$f(r)$};
		  \end{tikzpicture}
		  \caption{}\label{fig:f a}
    \end{subfigure}
    \hskip2em
    \begin{subfigure}{.45\linewidth}
    \centering
        \begin{tikzpicture}[xscale=0.6,yscale=0.6]
		    \draw [<-,thick] (0,3) -- (0,-3);
			\draw [->,thick] (-5,0) -- (3.5,0);
			\draw [-,thick] (0,0) to [out=65,in=180] (1.5,1.5) to [out=0,in=165] (3,0.25);
			
			\draw [-,thick] (-4.6,-0.5) to [out=-20,in=130] (-3.0,-2.0) to [out=-50, in=245] (0,0);
			
			\draw (1.5,-0.2) -- (1.5, 0.2);
			\draw (-2.25,-0.2) -- (-2.25, 0.2);
			\node [below] at (1.5,-0.2) {${r}^+$};
			\node [below] at (-2.0,-0.2) {${r}^-$};
			\node [right] at (3.5,0) {${r}$};
			\node [right] at (0,2.2) {$f'(r)$};
		  \end{tikzpicture}
		   \caption{}\label{fig:f b}
    \end{subfigure}
    \caption{(a) Potential function $f(r)$  for tensile force. $C^+$ and $C^-$ are two extreme values of $f$. (b) Cohesive tensile force.}\label{fig:ConvexConcave}
\end{figure}
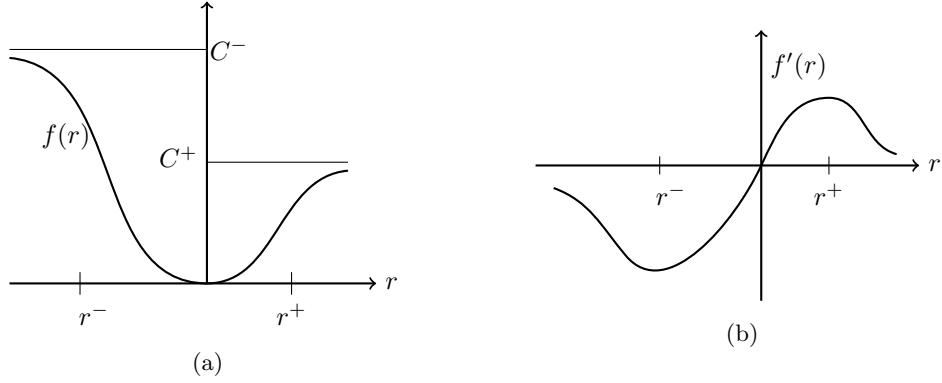

In this work we assume small deformation and work with linearized bond-strain. Let $D\subset \bbR^d$, for $d=2,3$, be the material domain. For a displacement field $\bu : D \times [0,T] \to \bbR^d$, the bond-strain between two material points $\bx, \by \in D$ is given by
\begin{align}\label{eq:bondstrain}
S(\by, \bx, t; \bu) = \frac{\bu(\by, t) - \bu(\bx, t)}{|\by - \bx|} \cdot \frac{\by - \bx}{|\by - \bx|}.
\end{align}
Let $\epsilon > 0$ be the size of horizon and $H_\epsilon(\bx) = \{ \by \in \bbR^d: |\by - \bx| < \epsilon\}$ be the neighborhood of a material point $\bx$. For pairwise interaction, we assume the following form of pairwise interaction potential
\begin{align}\label{eq:tensilepot}
\mathcal{W}^\epsilon(S(\by, \bx, t;\bu)) = \frac{J^\epsilon(|\by - \bx|)}{\epsilon|\by - \bx|} f(\sqrt{|\by - \bx|} S(\by, \bx, t;\bu)),
\end{align}
where $J^\epsilon(|\by - \bx|)$ is the influence function. We assume $J^\epsilon(|\by - \bx|) = J(|\by-\bx|/\epsilon)$ where $0\leq J(r) \leq M$ for $r<1$ and $J(r) = 0$ for $r\geq 1$. The potential $f$, see \autoref{fig:f a}, is assumed to be convex for small strains and becomes concave for larger strains. In the widely used prototypical micro-elastic brittle (PMB) peridynamic material, the strain vs force profile is linear up to some critical strain $S_c$ and is zero for any strain above $S_c$. In contrast, the peridynamic force given by $\partial_S \mathcal{W}^\epsilon$, is linear near zero strain and as the strain gets larger and reaches the critical strain, $S_c^+$ ($S_c^-$) for positive (negative) strain, the bond starts to soften, see \autoref{fig:f b}. For a given potential function $f$, the critical strain is given by $S_c^+ = \frac{r^+}{\sqrt{|\by - \bx|}}$ and $S_c^- = \frac{r^-}{\sqrt{|\by - \bx|}}$ where $r^+>0, r^-<0$ are the inflection points of the potential function $f$ as shown in \autoref{fig:f a}. 

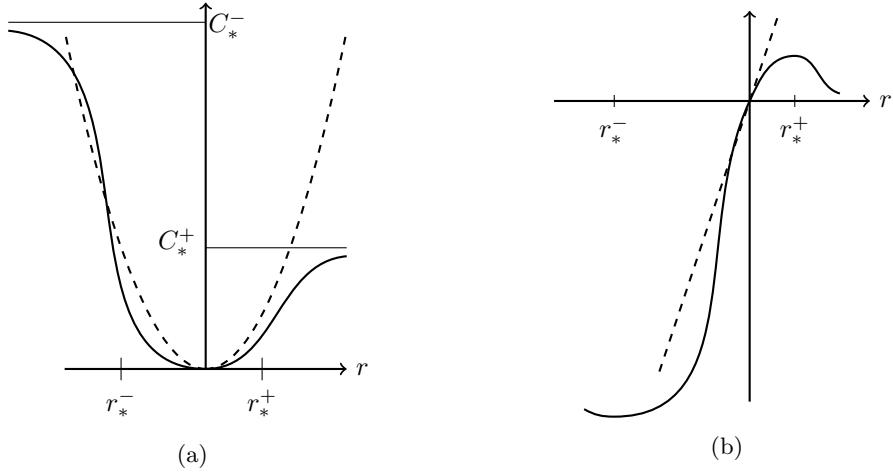
\begin{figure}
    \centering
    \begin{subfigure}{.45\linewidth}
    \centering
        \begin{tikzpicture}[xscale=0.75,yscale=0.75]
			\draw [<->,thick] (0,6.5) -- (0,0) -- (2.5,0);
			\draw [-,thick] (0,0) -- (-2.5,0);
			\draw [-] (0,2.15) -- (2.5,2.15);
			\draw [-] (0,6.15) -- (-3.5,6.15);
			\draw [-,thick] (0,0) to [out=0,in=-175] (2.5,2);
			\draw [-,thick] (-3.5,6) to [out=-5,in=180] (0,0);
			

			\draw [-,thick, dashed] (0,0) parabola (2.5,6);
			\draw [-,thick, dashed] (0,0) parabola (-2.5,6);
			
			\draw (1,-0.2) -- (1, 0.2);
			\node [below] at (1,-0.2) {${r}_\ast^+$};
			
			\draw (-1.5,-0.2) -- (-1.5, 0.2);
			\node [below] at (-1.5,-0.2) {${r}_\ast^-$};

			\node [right] at (2.5,0) {$r$};
			\node [left] at (0,2.250) {$C^+_\ast$};
			\node [left] at (0.9,6.1) {$C^-_\ast$};
		  \end{tikzpicture}
		  \caption{}\label{fig:g a}
    \end{subfigure}
    \hskip2em
    \begin{subfigure}{.45\linewidth}
    \centering
        \begin{tikzpicture}[xscale=0.4,yscale=0.4]
			\draw [<-,thick] (0,3) -- (0,-10);
			\draw [->,thick] (-6.5,0) -- (4,0);
			\draw [-,thick] (-5.5,-10.25) to [out=-30,in=180] (-4.5,-10.5) to [out=0, in=245 ] (0,0) to [out=65,in=180] (1.5,1.5)
			           to [out=0,in=165] (3,0.25);
			           
			\draw [-,thick,dashed] (-3,-9) -- (1,3);
			
			\draw (1.5,-0.2) -- (1.5, 0.2);
			\draw (-4.5,-0.2) -- (-4.5, 0.2);
			\node [below] at (1.5,-0.2) {${r}_\ast^+$};
			\node [below] at (-4.5,-0.2) {${r}_\ast^-$};
			\node [right] at (4,0) {${r}$};
		  \end{tikzpicture}
		  \caption{}\label{fig:g b}
    \end{subfigure}
    
    \caption{(a) Two types of potential function $g(r)$ for hydrostatic  force. Dashed line corresponds to the quadratic potential $g(r) = \beta r^2/2$. Solid line corresponds to the convex-concave type potential $g(r)$. For the convex-concave type potential, there are two special points $r^-_\ast$ and $r^+_\ast$ at which material points start to soften. $C^+_\ast$ and $C^-_\ast$ are two extreme values. (b) Hydrostatic forces.}
    \label{fig:ConvexConcaveFunctionG}
\end{figure}

The spherical or hydrostatic strain $\theta(\bx, t;\bu)$ at material point is given by 
\begin{align}\label{eq:sphericalstrain}
\theta(\bx,t;\bu)=\frac{1}{\epsilon^d \omega_d}\int_{H_\epsilon(\bx)} J^\epsilon(|\by-\bx|)S(\by,\bx,t;\bu){|\by-\bx|}\,d\by,
\end{align}
where $\omega_d$ is the volume of unit ball in dimension $d=2,3$. The potential for hydrostatic interaction is of the form
\begin{equation}\label{eq:hydropot}
\mathcal{V}^\epsilon(\theta(\bx,t;\bu))=\frac{g(\theta (\bx,t;\bu))}{\epsilon^2},
\end{equation}
where $g$ is the potential function associated to hydrostatic strain. Here $g$ can be of two types: 1) a quadratic function with only one well at zero strain, and 2) a convex-concave function with a  wells at the origin and at $\pm\infty$, see \autoref{fig:g a}. If $g$ is assumed to be quadratic then the force due to spherical strain is linear. If $g$ is a multi-well potential, the material softens as the hydrostatic strains exceeds critical value. For the convex-concave type $g$, the critical values are $0<\theta_c^+$ and $\theta^-_c< 0$ beyond which the force begins to soften is related to the inflection point $r^+_\ast$ and $r_\ast^-$ of $g$ as follows
\begin{equation}
\theta^+_c={r^+_\ast}, \qquad \theta^-_c={r_\ast^-}.
\label{eq:criticaltheta}
\end{equation}
The critical compressive hydrostatic strain where the force begins to soften for negative hydrostatic strain is chosen much larger in magnitude than $\theta_c^+$, i.e. $\theta^+_c << |\theta_c^-|$.

The finite element approximation has been applied to peridynamic fracture, however there remains a paucity of literature addressing the rigorous a-priori convergence rate of the finite element approximation to peridynamic problems in the presence of material failure. This aspect provides the motivation for the present work. In this paper we first prove existence of peridynamic evolutions taking values in $H^2(D;\bbR^d)\cap H^1_0(D;\bbR^d)$ that are twice differentiable in time, see \autoref{thm:existence}. We note that as these evolutions will become more fracture like as the region of nonlocal interaction decreases. These evolutions can be thought of as inner approximations to fracture evolutions. On passing to subsequences it is possible to show that the $H^2(D;\bbR^d)\cap H^1_0(D;\bbR^d)$ evolutions converge in the limit of vanishing non-locality to a limit solution  taking values in the space of special functions of bounded deformation SBD. Here the limit evolution has a well  defined Griffith fracture energy bounded by the initial data, see \cite{CMPer-Lipton} and \cite{Jha and Lipton 2019}. We show here that higher temporal regularity can be established if the body force changes smoothly in time. {Motivated by these considerations we develop finite element error estimates for solutions that take values in $H^2(D;\bbR^d)\cap H^1_0(D;\bbR^d)$ and for a bounded time interval.}

In this paper we obtain an a-priori $L^2$ error bound for the finite element {approximation of the displacement and velocity using a central in time discretization. } {Due to the nonlinear nature of the problem we get a convergence rate by using Lax Richtmyer stability together with consistency. Both stability and consistency are shown to follow from the Lipschitz continuity of the peridynamic force in $L^2(D;\bbR^d)$, see section \ref{consistent} and section  \ref{stable}.} The bound on the $L^2$ error is uniform in time and is given by $C_t \Delta t + C_s h^2/\epsilon^2$, where the constants  $C_t$ and $C_s$ are independent of $\Delta t$ and mesh size $h$, see \autoref{thm:convergence}. 
A more elaborate discussion of the a-priori bound is presented in section \autoref{s:convergence}.  {For the linearized model we obtain a stability condition on $\Delta t$, \autoref{thm:cflcondition}, that is of the same form as those given for linear local and nonlocal wave equations  \cite{CMPer-Karaa,CMPer-Guan}.} We demonstrate stability for the linearized model noting that for small strains the material behaves like a linear elastic material and that the stability of the linearized model is necessary for the stability of nonlinear model. 
{We believe a more constructive CFL stability condition is possible for the linear case and will pursue this in future work.}

{Previous work \cite{CMPer-JhaLiptonFD} treated spatially Lipschitz continuous solutions and addressed the finite difference approximation and obtained bounds on the $L^2$ error for {the displacement and velocity that are uniform in time and of the} form $C_t \Delta t + C_s h/\epsilon^2$, where constants the $C_t$ and $C_s$ are as before. For finite elements the convergence rate is seen to be slower than for the FEM model introduced here and is of order $h/\epsilon^2$ as opposed to $h^2/\epsilon^2$.  On the other hand the FEM method increases the computational work due to the inversion of the mass matrix.}

We carry out numerical experiments for dynamic crack propagation and obtain convergence rates for Plexiglass that are in line with the theory, see \autoref{s:numerical}. We also compare the Griffith's fracture energy with the peridynamic energy of the material softening zone we show good agreement between the two energies, see \autoref{sec: softzone}. Finite difference methods are less expensive than finite element approximations for nonlocal problems, however the latter offers more control on the accuracy of solution see, \cite{CMPer-Richard,CMPer-Littlewood,CMPer-Du6,CMPer-Gerstle, Du 2018b}. 

Here the a-priori $L^2$ convergence rates for the FEM given by \autoref{thm:convergence} include the effects of material degradation through the softening of material properties. The FEM simulations presented in this paper show that the material develops localized softening zones (region where bonds exceed critical tensile strain) as it deforms. This is in contrast to linear peridynamic models which are incapable of developing softening zones. For nonlinear peridynamic models with material failure the localization of zones of softening and damage is the hallmark of peridynamic modeling \cite{CMPer-Silling,States}, \cite{HaBobaru}, \cite{CMPer-Silling7}. {
One notes that the a-priori error involves $\epsilon$ in the denominator and in many cases $\epsilon$ is chosen small. However typical dynamic fracture experiments last only hundreds of microseconds and the a-priori error is controlled by the product of simulation time multiplied by $h^2/\epsilon^2$. So for material properties  characteristic of Plexiglass and  $\epsilon$  of size $4mm$,  the a-priori estimates predict a relative error of $1/10$ for simulations lasting around 100 microseconds.  
{We point out that the a-priori error estimates assume the appearance of nonlinearity anywhere in the computational domain. On the other hand  numerical simulation and independent theoretical estimates show that the nonlinearity concentrates along ``fat'' cracks of finite length  and width equal to $\epsilon$, see \cite{CMPer-Lipton,CMPer-Lipton3}. Moreover the remainder of the computational domain is seen to behave linearly and to leading order can be modeled as a linear elastic material up to an error proportional to $\epsilon$, see [Proposition 6, \cite{Jha and LiptonL}]. Future work will use these observations to focus on adaptive implementation and a-posteriori estimates.}
} A-posteriori  convergence for FEM models of peridynamics with material degradation  can be seen in the work \cite{CMPer-Richard}, \cite{CMPer-Chen}, \\\cite{CMPer-Ren}. For other nonlinear and nonlocal models adaptive  mesh refinement within FE framework for nonlocal models has been explored in \cite{CMPer-Du6} and convergence of the adaptive FE approximation is rigorously shown. A-posteriori error analysis of linear nonlocal models is carried out in \cite{CMPer-Du7}. 

The paper is organized as follows. We introduce the equation of motion in  \autoref{s:nonlocal dynamics} and present the Lipschitz continuity of the force,  existence of peridynamic solution, and the higher temporal regularity necessary for the finite element error analysis. 
In \autoref{s:fem}, we consider the finite element discretization. 
We prove the stability of a semi-discrete approximation in \autoref{ss:semifem}. In \autoref{s:fullfem}, we analyze the fully discrete approximation and obtain an a-priori bound on errors. The stability of the fully discrete approximation linearized peridynamic force is shown in \autoref{ss:fullfemstability}. We present our numerical experiments in \autoref{s:numerical}. 
Proofs of the Lipschitz bound on the peridynamic force and higher temporal regularity of solutions is provided in \autoref{s:proofs}. 
In \autoref{s:conclusions} we present our conclusions.

We conclude the introduction by listing the notation used throughout the paper.
We denote material domain as $D$ where $D\subset \bbR^d$ for $d=2,3$. Points and vectors in $\bbR^d$ are denoted as bold letters. Some of the key notations are as follows

{\vskip 2mm}

\begin{tabular}{lp{0.65\textwidth}}
  $[0,T]$ & Time domain \\
  $\epsilon$ & Size of horizon \\
  $\rho$ & Density \\
  $H_\epsilon(\bx)$ & Horizon of $\bx \in D$, a ball of radius $\epsilon$ centered at $\bx$  \\
  $\omega_d$ & Volume of unit ball in dimension $d=2,3$ \\
  $\omega(\bx) \in [0,1]$ & Boundary function defined on $D$ taking value $1$ in the interior and smoothly decaying to $0$ as $\bx$ approaches $\partial D$ \\
  $\bu$ & Displacement field defined over $D \times [0,T]$. We may also use notation $\bu$ to denote field defined over just $D$\\
  $\bu_0, \bv_0$ & Initial condition on displacement \\
  $\bb$ & Body force defined over $D\times [0,T]$ \\
  $\be_{\by-\bx}$ & The unit vector pointing from a point $\by$ to the point $\bx$\\
  $S = S(\by, \bx, t; \bu)$ & Bond strain $S=\frac{\bu(\by,t)-\bu(\bx,t)}{|\by-\bx|}\cdot\be_{\by-\bx}$. We may also use $S(\by,\bx;\bu)$ if $\bu$ is a filed defined over just $D$ \\
  $\theta = \theta(\bx,t;\bu)$ & Spherical or hydrostatic strain. We may also use $\theta(\bx;\bu)$ if $\bu$ is a filed defined over just $D$ \\
  $S^+_c, S^-_c$ & Critical bond strain\\
  $\theta^+_c, \theta^-_c$ & Critical hydrostatic strain\\
  $J^\epsilon(r) = J(r/\epsilon)$ & Influence function where $J$ is integrable with $J(r) = 0$ for $r\geq 1$ and $0\leq J(r) \leq M$ for $r < 1$\\
  $\bar{J}_\alpha$ & Moment of function $J$ over $H_1(\bzero)$ with weight $1/(\omega_d |\bxi|^\alpha)$ \\
  $f, g$ & Potential functions for pairwise and state-based interaction \\
  $\mathcal{W}^\epsilon, \mathcal{V}^\epsilon$ & Pairwise and state-based potential energy density \\
  $PD^\epsilon(\bu(t))$ & Total peridynamic potential energy at time $t$ \\
  $\mathcal{E}^\epsilon(\bu)(t)$ & Total dynamic energy at time $t$ \\
  $\mathcal{L}^\epsilon, \mathcal{L}^\epsilon_T, \mathcal{L}^\epsilon_D$ & total peridynamic force, pairwise peridynamic force, and state-based peridynamic force respectively \\
  $a^\epsilon(\bu, \bv)$ & Nonlinear operator where $\bu, \bv$ are vector fields over $D$ \\
  $a^\epsilon_T, a^\epsilon_D$ & Nonlinear pairwise and state-based operator\\
  $||\cdot ||,||\cdot ||_\infty, ||\cdot ||_n$ & $L^2$ norm over $D$, $L^\infty$ norm over $D$, and Sobolev $H^n$ norm over $D$ (for $n=1,2$) respectively \\
  $h, \Delta t$ & Size of mesh and size of time step\\
  $\mathcal{T}_h$ & Triangulation of $D$ given by triangular/tetrahedral elements \\
  $\mathcal{I}_h$ & Continuous piecewise linear interpolation operator on $\mathcal{T}_h$ \\
  $W $ & Space of functions in $H^2(D;\bbR^d)$ such that trace of function is zero on boundary $\partial D$, i.e. $W = H^2(D;\bbR^d) \cap H^1_0(D;\bbR^d)$\\
  $V_h$ & Space of continuous piecewise linear interpolations on $\mathcal{T}_h$ \\
  $\phi_i$ & Interpolation function of mesh node $i$ \\
  $\br_h(\bu)$ & Finite element projection of $\bu$ onto $V_h$\\
  $E^k$ & Total error in mean square norm at time step $k$\\
  $\bu^k_h, \bv^k_h$ & Approximate displacement and velocity field at time step $k$\\
  $\bu^k, \bv^k$ & Exact displacement and velocity field at time step $k$\\
\end{tabular}\\


\section{Equation of motion, existence, uniqueness, and higher regularity}\label{s:nonlocal dynamics}
We assume $D$ to be an open set with $C^1$ boundary.  To enforce zero displacement boundary conditions at $\partial D$ and to insure a well posed evolution we introduce the boundary function $\omega(\bx)$. This function is introduced as a factor into the potentials $\mathcal{W}^\epsilon$ and $\mathcal{V}^\epsilon$. Here the boundary function takes value $1$ in the interior of domain and is zero on the boundary.  We assume $\sup_{\bx} |\nabla \omega(\bx)| < \infty$ and $\sup_{\bx} |\nabla^2 \omega(\bx)| < \infty$ in our analysis. The hydrostatic strain is modified to include the boundary and is given by
\begin{align}\label{eq:sphericalstrainmod}
\theta(\bx,t;\bu)=\frac{1}{\epsilon^d \omega_d}\int_{H_\epsilon(\bx)} \omega(\by) J^\epsilon(|\by-\bx|)S(\by,\bx,t;\bu){|\by-\bx|}\,d\by.
\end{align}
The peridynamic potentials \autoref{eq:tensilepot} and \autoref{eq:hydropot} are modified to {\em see the boundary} follows
\begin{align}
\mathcal{W}^\epsilon(S(\by, \bx, t;\bu)) &= \omega(\bx) \omega(\by) \frac{J^\epsilon(|\by - \bx|)}{\epsilon|\by - \bx|} f(\sqrt{|\by - \bx|} S(\by, \bx, t;\bu)), \label{eq:tensilepotmod} \\
\mathcal{V}^\epsilon(\theta(\bx,t;\bu)) &= \omega(\bx) \frac{g(\theta (\bx,t;\bu))}{\epsilon^2} \label{eq:hydropotmod}.
\end{align}
We assume that potential function $f$ is at least 4 times differentiable and satisfies following regularity condition:
\begin{align}
C^f_0 := \sup_r |f(r)| < \infty, \qquad C^f_i := \sup_r |f^{(i)}(r)| < \infty, \quad \forall i = 1,2,3,4. 
\end{align}
If potential function $g$ is convex-concave type then we assume that $g$ satisfies same regularity condition as $f$. We denote constants $C^g_i$, for $i=0,1,..,4$, similar to $C^f_i$ above.

The total potential energy at time $t$ is given by
\begin{equation}\label{eq:totalenergy}
\begin{aligned}
PD^\epsilon(\bu(t))=\frac{1}{\epsilon^d \omega_d}\int_D \int_{H_\epsilon(\bx)} |\by-\bx|\mathcal{W}^\epsilon(S(\by,\bx,t;\bu))\,d\by d\bx\\
+\int_D \mathcal{V}^\epsilon(\theta(\bx,t;\bu))\,d\bx,
\end{aligned}
\end{equation}
where potential $\mathcal{W}^\epsilon$ and $\mathcal{V}^\epsilon$ are described above. The material is assumed to be homogeneous and the density is given by $\rho$. The applied body force is denoted by $\bb(\bx,t)$. We define the {Lagrangian}  $${\rm{L}}(\bu,\partial_t \bu,t)=\frac{\rho}{2}||\dot \bu||^2 - PD^\epsilon(\bu(t))+\int_D \bb(t)\cdot \bu(t) d\bx,$$
here $\dot \bu$ is the velocity given by the time derivative of $\bu$. Applying the {principal of least action} together with a straight forward calculation (see for example \cite{CMPer-Lipton5} for detailed derivation) gives the nonlocal dynamics 
\begin{equation}\label{eq:equationofmotion}
\begin{aligned}
\rho \ddot{\bu}(\bx,t)=\mathcal{L}^\epsilon(\bu)(\bx,t)+\bb(\bx,t),\hbox{  for  $\bx\in D$},
\end{aligned}
\end{equation}
where 
\begin{align}\label{eq:totalforce}
\mathcal{L}^\epsilon(\bu)(\bx,t) = \mathcal{L}^\epsilon_T(\bu)(\bx,t) + \mathcal{L}^\epsilon_D(\bu)(\bx,t),
\end{align}
$\mathcal{L}^\epsilon_T(\bu)$ is the peridynamic force due to bond-based interaction and is given by
\begin{align}\label{eq:nonlocforcetensile}
&\mathcal{L}^\epsilon_T(\bu)(\bx,t) \notag \\
&=\frac{2}{\epsilon^d \omega_d}\int_{H_\epsilon(\bx)} \omega(\bx) \omega(\by) \frac{J^\epsilon(|\by-\bx|)}{\epsilon|\by-\bx|}\partial_S f(\sqrt{|\by-\bx|}S(\by,\bx,t;\bu))\be_{\by-\bx}\,d\by,
\end{align}
and $\mathcal{L}^\epsilon_D(\bu)$ is the peridynamic force due to state-based interaction and is given by
\begin{align}\label{eq:nonlocforcedevia}
&\mathcal{L}^\epsilon_D(\bu)(\bx,t)\notag \\
&=\frac{1}{\epsilon^d \omega_d}\int_{H_\epsilon(\bx)} \omega(\bx) \omega(\by) \frac{J^\epsilon(|\by-\bx|)}{\epsilon^2}\left[\partial_\theta g(\theta(\by,t;\bu))+\partial_\theta g(\theta(\bx,t;\bu))\right]\be_{\by-\bx}\,d\by.
\end{align}

The dynamics is complemented with the initial data 
\begin{align}\label{eq:idata}
\bu(\bx,0)=\bu_0(\bx), \qquad  \partial_t \bu(\bx,0)=\bv_0(\bx).
\end{align}
We prescribe zero Dirichlet boundary condition on the boundary $\partial D$
\begin{align}\label{eq:bc}
\bu(\bx) = \bzero \qquad \forall \bx \in \partial D.
\end{align}
We extend the zero boundary condition outside $D$ to whole $\bbR^d$. In our analysis we will assume mass density $\rho = 1$ without loss of generality. 

\subsection{Existence of solutions and higher regularity in time}\label{s:existence}
We recall that the space $H^n_0(D;\bbR^d)$ is the closure in the $H^n$ norm of the  functions that are infinitely differentiable with compact support in $D$. For suitable initial conditions and body force we show that solutions exist in 
\begin{equation}
\label{W}
W=H^2(D;\bbR^d) \cap H^1_0(D;\bbR^d)=\{v\in H^2(D;\bbR^d):\gamma v =0,\hbox{  on  }\partial D\}
\end{equation}
where $\gamma$ is the trace of the function $v$ on the boundary of $D$. 
We will assume that $\bu\in W$ is extended by zero outside $D$. We first exhibit the Lipschitz continuity property and boundedness of the peridynamic force for displacements in $W$. We will then apply [Theorem 3.2, \cite{CMPer-JhaLipton3}] to conclude the existence of unique solutions.

We note the following Sobolev embedding properties of $H^2(D;\bbR^d)$ when $D$ is a $C^1$ domain. 
\begin{itemize}
\item From Theorem 2.72 of \cite{MAFa-Demengel}, there exists a constant $C_{e_1}$ independent of $\bu \in H^2(D;\bbR^d)$ such that 
\begin{align}
||\bu||_{\infty} \leq C_{e_1} ||\bu||_{2} . \label{eq:sob embedd 1}
\end{align}

\item Further application of standard embedding theorems \\(e.g., Theorem 2.72 of \cite{MAFa-Demengel}), shows there exists a constant $C_{e_2}$ independent of $\bu$ such that
\begin{align}
||\nabla \bu||_{L^q(D;\bbR^{d\times d})} \leq C_{e_2} ||\nabla \bu||_{1} \leq C_{e_2} ||\bu||_{2}, \label{eq:sob embedd 2}
\end{align}
for any $q$ such that $2\leq q< \infty$ when $d=2$ and $2\leq q \leq 6$ when $d=3$.
\end{itemize}

We have following result which show the Lipschitz continuity property of a peridynamic force $\mathcal{L}^\epsilon$. 
{\vskip 2mm}
\begin{theorem}\label{thm:lipschitzproperty}
\textbf{Lipschitz continuity of peridynamic force}\\
Let $f$ be a convex-concave function satisfying $C^f_i <\infty$ for $i=0,..,4$ and let $g$ either be a quadratic function, or $g$ be a convex-concave function with $C^g_i < \infty$ for $i=0,..,4$. Also, let boundary function $\omega : D\to [0,1]$ be such that $\sup_{\bx\in D} |\nabla \omega(\bx)| <\infty$ and $\sup_{\bx\in D} |\nabla^2 \omega(\bx)| <\infty$. Then, for any $\bu,\bv\in W$, we have
\begin{align}
||\mathcal{L}^\epsilon(\bu) - \mathcal{L}^\epsilon(\bv)||_2 &\leq \dfrac{\bar{L}_1 (1 + ||\bu||_2 + ||\bv||_2)^2}{\epsilon^3} ||\bu - \bv||_2,
\end{align}
where constant $\bar{L}_1$ does not depend on $\epsilon$ and $\bu,\bv$. Also, for $\bu \in W$, we have
\begin{align}
||\mathcal{L}^\epsilon(\bu)||_2 &\leq \dfrac{\bar{L}_2 (||\bu||_2 + ||\bu||_2^2)}{\epsilon^{5/2}},
\end{align}
where constant $\bar{L}_2$ does not depend on $\epsilon$ and $\bu$. 
\end{theorem}
{\vskip 2mm}

Now let $T>0$ be any positive number, a straight-forward application of [Theorem 3.2, \cite{CMPer-JhaLipton3}] gives:
{\vskip 2mm}
\begin{theorem}\label{thm:existence} 
\textbf{Existence and uniqueness of solutions over finite time intervals}\\
Let $f$, $g$, and $\omega$ satisfy the hypothesis of \autoref{thm:lipschitzproperty}. For any initial condition $\bu_0,\bv_0 \in W$, time interval $I_0=(-T,T)$, and right hand side $\bb(t)$  continuous in time for $t\in I_0$ such that $\bb(t)$ satisfies $\sup_{t\in I_0} ||\bb(t)||_2<\infty$, there is a unique solution $\bu(t)\in C^2(I_0;W)$ of peridynamic equation \ref{eq:equationofmotion}. Also, $\bu(t)$ and $\dot{\bu}(t)$ are Lipschitz continuous in time for $t\in I_0$.
\end{theorem}
{\vskip 2mm}

We can also show higher regularity in time of evolutions under suitable assumptions on the body force:
{\vskip 2mm}
\begin{theorem}\label{thm:higherregularity} 
\textbf{Higher regularity}\\
Suppose the initial data and righthand side $\bb(t)$ satisfy the hypothesis of \autoref{thm:existence} and suppose further that  $\dot{\bb}(t)$ exists and is continuous in time for $t\in I_0$ and \,$\sup_{t\in I_0} ||\dot{\bb}(t)||_2 < \infty$. Then 
$\bu \in C^3(I_0; W)$ and 
\begin{align}
|| \partial^3_{ttt} \bu(\bx,t)||_2 &\leq  \dfrac{C(1 + \sup_{s \in I_0} ||\bu(s)||_2)^2}{\epsilon^3} \sup_{s\in I_0} ||\partial_t \bu(s)||_2 + ||\dot{\bb}(\bx,t)||_2,
\end{align}
where $C$ is a positive constant independent of $\bu$.

\end{theorem}
{\vskip 2mm}

The proofs of Theorems \ref{thm:lipschitzproperty}, and \ref{thm:higherregularity} are given in \autoref{s:proofs}. For future reference, we note that for any $\bu,\bv \in L^2_0(D;\bbR^d)$, we have
\begin{align}\label{eq:lipschitzproperty l2}
||\mathcal{L}^\epsilon(\bu) - \mathcal{L}^\epsilon(\bv) || &\leq \frac{L}{\epsilon^2} ||\bu - \bv||.
\end{align}
where constant $L$ is given by
\begin{align}\label{eq:lipschitzconstant l2}
L&:= \begin{cases}
4(C^f_2 \bar{J}_1 + C^g_2\bar{J}_0^2) \qquad \text{if $g$ is a convex-concave type}, \\
4(C^f_2 \bar{J}_1 + g''(0)\bar{J}_0^2) \qquad \text{if $g$ is a quadratic function},
\end{cases}
\end{align}
and $\bar{J}^\alpha = (1/\omega_d)\int_{H_1(\bzero)} J(|\bxi|)/|\bxi|^\alpha d\bxi$ . 

\subsection{Weak form}\label{ss:weakform}
We multiply \autoref{eq:equationofmotion} by a test function $\butilde$ in $H^1_0(D;\bbR^d)$ and integrate over $D$ to get
\begin{align}
(\ddot{\bu}(t),\butilde) = (\mathcal{L}^\epsilon(\bu(t)), \butilde) + (\bb(t),\butilde).
\end{align}

We have the following integration by parts formula:
\begin{lemma}\label{lem:operatora}
For any $\bu, \bv \in L^2_0(D;\bbR^d)$ we have
\begin{align}
(\mathcal{L}^\epsilon(\bu), \bv)  = - a^\epsilon(\bu, \bv), \label{eq:bypartsformula}
\end{align}
where
\begin{align}\label{eq:operatora decomp}
a^\epsilon(\bu, \bv) = a^\epsilon_T(\bu, \bv) + a^\epsilon_D(\bu, \bv)
\end{align}
and
\begin{align}\label{eq:operatora TandD}
a^\epsilon_T(\bu, \bv) &= \dfrac{1}{\epsilon^{d+1}\omega_d} \int_D \int_D \omega(\bx) \omega(\by) J^\epsilon(|\by - \bx|) \notag \\
&\qquad \qquad\partial_S f(\sqrt{|\by - \bx|} S(\by,\bx;\bu)) S(\by, \bx;\bv) d\by d\bx, \notag \\
a^\epsilon_D(\bu, \bv) &= \dfrac{1}{\epsilon^2} \int_D \omega(\bx) g'(\theta(\bx;\bu)) \theta(\bx;\bv) d\bx.
\end{align}
\end{lemma}
The proof of above lemma is identical to the proof of Lemma 4.2 in \cite{CMPer-Lipton5}.

Using the above Lemma, the weak form of the peridynamic evolution is given by 
\begin{align}\label{eq:weakform}
(\ddot{\bu}(t),\butilde) + a^\epsilon(\bu(t), \butilde) = (\bb(t),\butilde).
\end{align}

\noindent\textit{Total dynamic energy: }We define the total dynamic energy as follows
\begin{align}\label{eq:totaldynamicenergy}
\mathcal{E}^\epsilon(\bu)(t) = \frac{1}{2} ||\dot{\bu}(t)||^2_{L^2} + PD^\epsilon(\bu(t)),
\end{align}
where $PD^\epsilon$ is defined in \autoref{eq:totalenergy}. Time derivative of total energy satisfies
\begin{align}\label{eq:energyidentity}
\frac{d}{d t} \mathcal{E}^\epsilon(\bu)(t) = (\ddot{\bu}(t), \dot{\bu}(t)) + a^\epsilon(\bu(t), \dot{\bu}(t)).
\end{align}

\textbf{Remark.} It is readily verified that the peridynamic force and energy are bounded for all functions in $L^2(D;\bbR^d)$. Here the bound on the force follows from the Lipschitz property of the force in $L^2(D;\bbR^d)$, see, \autoref{eq:lipschitzproperty l2}.  The peridynamic force is also bounded for functions $\bu$ in $H^1(D;\bbR^d)$. This again follows  from the Lipschitz property of the force in $H^1(D;\bbR^d)$ using arguments established in \autoref{s:proofs}.  The boundedness of the energy $PD^\epsilon(\bu)$  in both $L^2(D;\bbR^d)$ and $H^1(D;\bbR^d)$ follows from the boundedness of the bond potential  energy $\mathcal{W}^\epsilon(S(\by, \bx, t; \bu))$ and $\mathcal{V}^\epsilon(\theta(\bx, t; \bu))$ used in the definition of $PD^\epsilon(\bu)$, see \autoref{eq:tensilepotmod} and \autoref{eq:hydropotmod}. More generally this also shows that $PD^\epsilon(\bu)<\infty$ for $\bu\in L^1(D;\bbR^d)$.

We next discuss the spatial and the time discretization of peridynamic equation. 

\section{Finite element approximation}\label{s:fem}
Let $V_h$ be given by linear continuous interpolations
over  tetrahedral or triangular elements $\mathcal{T}_h$ where $h$ denotes the size of finite element mesh. Here we  assume the elements are conforming and the finite element mesh is shape regular and $V_h\subset H^1_0(D;\bbR^d)$.

For a continuous function $\bu$ on $\bar{D}$, $\mathcal{I}_h(\bu)$ is the continuous piecewise linear interpolant on $\mathcal{T}_h$. It is given by 
\begin{align}
\mathcal{I}_h(\bu)\bigg\vert_{T} = \mathcal{I}_T(\bu) \qquad \forall T\in \mathcal{T}_h,
\end{align}
where $\mathcal{I}_T(\bu)$ is the local interpolant defined over finite element $T$ and is given by
\begin{align}
\mathcal{I}_T(\bu) = \sum_{i=1}^n \bu(\bx_i)\phi_i.
\end{align}
Here $n$ is the number of vertices in an element $T$, $\bx_i$ is the position of vertex $i$, and $\phi_i$ is the linear interpolant associated to vertex $i$. 

Application of Theorem 4.4.20 and remark 4.4.27 in \cite{MANa-Susanne} gives 
\begin{align}\label{eq:interpolationerror}
&&|| \bu - \mathcal{I}_h(\bu) || \leq c h^2 || \bu ||_2, \hbox{         } \qquad \forall \bu \in W.
\end{align}

Let $\br_h(\bu)$ denote the projection of $\bu\in W$ on $V_h$. For the $L^2$ norm it is defined as
\begin{align}
||\bu - \br_h(\bu)|| &= \inf_{\tilde{\bu}\in V_h} ||\bu - \tilde{\bu}|| . \label{eq:projection}
\end{align}
and satisfies 
\begin{align}\label{eq:projectionorthogonal}
(\br_h(\bu), \tilde{\bu}) = (\bu, \tilde{\bu}), \qquad \forall \tilde{\bu} \in V_h.
\end{align}

Since $\mathcal{I}_h(\bu) \in V_h$, and \autoref{eq:interpolationerror} we see that
\begin{align}
&&||\bu - \br_h(\bu)|| \leq c h^2 || \bu ||_2,\hbox{     } \qquad \forall \bu \in W \label{eq:projectionerror}
\end{align}

\subsection{Semi-discrete approximation}\label{ss:semifem}
Let $\bu_h(t) \in V_h$ be the approximation of $\bu(t)$ satisfying following for all $t\in [0,T]$
\begin{align}
(\ddot{\bu}_h, \tilde{\bu} ) + a^\epsilon(\bu_h(t), \tilde{\bu}) &= ( \bb(t), \tilde{\bu} ), \qquad \forall \tilde{\bu} \in V_h\label{eq:feweak}.
\end{align}
We have following result:
{\vskip 2mm}
\begin{theorem}\label{thm:stabilitysemi}
\textbf{Energy stability of semi-discrete approximation}\\
The semi-discrete scheme is stable and the energy $\mathcal{E}^\epsilon(\bu_h)(t)$, defined in \autoref{eq:totaldynamicenergy}, satisfies the following bound
\begin{align*}
\mathcal{E}^\epsilon(\bu_h)(t) &\leq \left[ \sqrt{\mathcal{E}^\epsilon(\bu_h)(0)} + \int_0^t ||\bb(\tau)|| d\tau \right]^2.
\end{align*}
\end{theorem}
{\vskip 2mm}
We note that while proving the stability of semi-discrete scheme corresponding to nonlinear peridynamics, we do not require any assumption on the strain $S(\by,\bx, t;\bu_h)$. The proof is similar to [Section 6.2, \cite{CMPer-Lipton}]. 

\begin{proof}
Letting $\tilde{\bu} = \dot{\bu}_h(t)$ in \autoref{eq:feweak} and noting the identity \autoref{eq:energyidentity}, we get
\begin{align}\label{eq:ineq stab 1}
\dfrac{d}{dt} \mathcal{E}^\epsilon(\bu_h)(t) = (\bb(t), \dot{\bu}_h(t)) \leq ||\bb(t)||\, ||\dot{\bu}_h(t)||.
\end{align}
We also have 
\begin{align*}
||\dot{\bu}_h(t)|| \leq 2 \sqrt{\frac{1}{2} ||\dot{\bu}_h||^2 + PD^\epsilon(\bu_h(t))} = 2 \sqrt{\mathcal{E}^\epsilon(\bu_h)(t)}
\end{align*}
where we used the fact that $PD^\epsilon(\bu)(t)$ is nonnegative. We substitute above inequality in \autoref{eq:ineq stab 1} to get
\begin{align*}
\dfrac{d}{dt} \mathcal{E}^\epsilon(\bu_h)(t) &\leq 2 \sqrt{\mathcal{E}^\epsilon(\bu_h)(t)} \, ||\bb(t)||.
\end{align*}
We fix $\delta > 0$ and define $A(t)$ as $A(t) = \mathcal{E}^\epsilon(\bu_h)(t) + \delta$. Then, from the above equation we easily have
\begin{align*}
\dfrac{d}{dt} A(t) &\leq 2 \sqrt{A(t)} \, ||b(t)|| \quad \Rightarrow \dfrac{1}{2} \dfrac{\frac{d}{dt} A(t)}{\sqrt{A(t)}} \leq ||\bb(t)||.
\end{align*}
Noting that $\frac{1}{\sqrt{a(t)}}\frac{d a(t)}{dt}= 2\frac{d}{dt} \sqrt{a(t)} $, integrating from $t=0$ to $\tau$ and relabeling $\tau$ as $t$, we get
\begin{align*}
\sqrt{A(t)} &\leq \sqrt{A(0)} + \int_0^t ||\bb(s)|| ds.
\end{align*}
Proof is complete once we let $\delta \to 0$ and take the square of both sides.
\end{proof}

\section{Central difference time discretization}\label{s:fullfem}
In \autoref{s:convergence} we calculate the convergence rate for the central difference time discretization of the fully nonlinear problem. We then present a CFL like condition on the time step $\Delta t$ for the linearized peridynamic equation in \autoref{ss:fullfemstability}.

At time step $k$, the exact solution is given by $(\bu^k, \bv^k)$ where $\bv^k = \partial \bu^k/\partial t$, and their projection onto $V_h$ is given by $(\br_h(\bu^k), \br_h(\bv^k))$. The solution of fully discrete problem at time step $k$ is given by $(\bu^k_h, \bv^k_h)$. 

We approximate the initial data on displacement $\bu_0$ and velocity $\bv_0$ by their projections $\br_h(\bu_0)$ and $\br_h(\bv_0)$. Let $\bu^0_h = \br_h(\bu_0)$ and $\bv^0_h = \br_h(\bv_0)$. For $k\geq 1$, $(\bu^k_h, \bv^k_h)$ satisfies, for all $\tilde{\bu} \in V_h$,
\begin{align}
\left( \dfrac{\bu^{k+1}_h - \bu^k_h}{\Delta t}, \tilde{\bu} \right) &= (\bv^{k+1}_h, \tilde{\bu}), \notag \\
\left( \dfrac{\bv^{k+1}_h - \bv^k_h}{\Delta t}, \tilde{\bu} \right) &= (\mathcal{L}^\epsilon(\bu^k_h), \tilde{\bu} ) + (\bb^k_h, \tilde{\bu}) \label{eq:forward},
\end{align}
where we have denoted the projection of $\bb(t^k)$, i.e. $\br_h(\bb(t^k))$, as $\bb^k_h$. Combining the two equations delivers central difference equation for $\bu^k_h$. We have
\begin{align}\label{eq:central}
\left( \dfrac{\bu^{k+1}_h - 2 \bu^k_h + \bu^{k-1}_h}{\Delta t^2}, \tilde{\bu} \right) &= (\mathcal{L}^\epsilon(\bu^k_h), \tilde{\bu} ) + (\bb^k_h, \tilde{\bu}), \qquad \forall \tilde{\bu} \in V_h.
\end{align}
For $k=0$, we have $\forall \tilde{\bu} \in V_h$
\begin{align}\label{eq:centralzero}
\left( \dfrac{\bu^1_h - \bu^0_h}{\Delta t^2}, \tilde{\bu}\right) &= \dfrac{1}{2}(\mathcal{L}^\epsilon(\bu^0_h), \tilde{\bu})+ \dfrac{1}{\Delta t} (\bv^0_h, \tilde{\bu}) + \dfrac{1}{2}(\bb^0_h, \tilde{\bu}).
\end{align}


\subsection{Implementation details}
For completeness we describe the implementation of the time stepping method using FEM interpolants.  Let $\bolds{N}$ be the shape tensor then $\bu^k_h, \tilde{\bu} \in V_h$ are given by
\begin{align}
\bu^k_h = \bolds{N} \bolds{U}^k, \qquad \tilde{\bu}  = \bolds{N} \tilde{\bolds{U}},
\end{align}
where $\bolds{U}^k$ and $\tilde{\bolds{U}}$ are $ Nd$ dimensional vectors, where $N$ is the number of nodal points in the mesh and $d$ is the dimension.  

From \autoref{eq:central}, for all $\tilde{\bolds{U}} \in \bbR^{Nd}$ with elements of $\tilde{\bolds{U}}$ zero on the boundary, then the following holds for $k\geq 1$
\begin{align}\label{eq:feweakequation}
\left[\bolds{M}\frac{\bolds{U}^{k+1} - 2\bolds{U}^{k} + \bolds{U}^{k-1}}{\Delta t^2} \right] \cdot \tilde{\bolds{U}} = \bolds{F}^k \cdot \tilde{\bolds{U}} .
\end{align}
Here the mass matrix $\bolds{M}$ and force vector $\bolds{F}^k$ are given by
\begin{align}\label{eq:MandF}
\bolds{M} &:= \int_D \bolds{N}^T \bolds{N} d\bx,\notag \\
\bolds{F}^k &:= \bolds{F}^k_{pd}+ \int_D \bolds{N}^T \bb(\bx,t^k) d\bx,
\end{align}
where $\bolds{F}^k_{pd}$ is defined by
\begin{align}\label{eq:Fpd}
\bolds{F}^k_{pd} &:= \int_D \bolds{N}^T (\mathcal{L}^\epsilon(\bu^k_h)(\bx)) d\bx.
\end{align}
We remark that a similar equation holds for $k=0$.

At the time step $k$ we must invert ${\bolds{M}}$ to solve for ${\bolds{U}}^{k+1}$ using
\begin{align}\label{eq:weakformsolve}
{\bolds{U}}^{k+1} = \Delta t^2 {\bolds{M}}^{-1} {\bolds{F}}^k + 2 \bolds{U}^k - {\bolds{U}}^{k-1}.
\end{align}
As is well known this inversion amounts to an increase of computational complexity associated with discrete approximation of the weak formulation of the evolution. Further, the matrix-vector multiplication $\bolds{M}^{-1} \bolds{F}^k$ needs to be carried out at each time step. On the other hand the quadrature error in the computation of the force vector $\bolds{F}^k_{pd}$ is reduced when using the weak form. 

We next show the convergence of approximation.

\subsection{Convergence of approximation}\label{s:convergence}
In this section, we prove the uniform bound on the error and show that the approximate solution converges to the exact solution with rate given by $C_t \Delta t + C_s h^2/\epsilon^2$. Here horizon $\epsilon > 0$ is assumed to be fixed. We first compare the exact solution with its projection in $V_h$ and then compare the projection with the approximate solution. We further divide the calculation of error between the projection and the approximate solution in two parts, namely consistency analysis and error analysis.

Error $E^k$ is given by
\begin{align*}
E^k := ||\bu^k_h - \bu(t^k)|| + ||\bv^k_h - \bv(t^k)||.
\end{align*}
The error is split into two parts as follows
\begin{align*}
E^k &\leq \left( ||\bu^k - \br_h(\bu^k)|| + ||\bv^k - \br_h(\bv^k)|| \right) + \left( || \br_h(\bu^k) - \bu^k_h|| + ||\br_h(\bv^k) - \bv^k_h|| \right),
\end{align*}
where the first term is the error between the exact solution and projection, and the second term is the error between the projection and approximate solution. Let 
\begin{align}\label{eq:ek}
\be^k_h(\bu) := \br_h(\bu^k) - \bu^k_h \quad \text{ and  } \quad \be^k_h(\bv) := \br_h(v^k) - \bv^k_h
\end{align}
and
\begin{align}\label{eq:eknorm}
e^k := ||\be^k_h(\bu)|| + ||\be^k_h(\bv)||.
\end{align}

Using \autoref{eq:projectionerror}, we have
\begin{align}\label{eq:Ekineq}
E^k &\leq C_p h^2 + e^k,
\end{align}
where
\begin{align}\label{eq:constCp}
C_p :=  c\left[ \sup_t ||\bu(t)||_2 + \sup_t ||\dfrac{\partial \bu(t)}{\partial t}||_2 \right].
\end{align}

We have the following a-priori convergence rate given by
{\vskip 2mm}
\begin{theorem}\label{thm:convergence}
\textbf{Convergence of Central difference approximation}\\
Let $(\bu,\bv)$ be the exact solution of the peridynamic equation \ref{eq:equationofmotion}. Let $(\bu^k_h, \bv^k_h)$ be the FE solution of \autoref{eq:forward}. If $\bu,\bv \in C^2([0,T]; W)$, then the scheme is consistent and the error $E^k$ satisfies following bound
\begin{align}\label{eq:Ekbound}
&\sup_{k \leq T/\Delta t} E^k \notag \\
&= C_p h^2 + \exp[T(1+L/\epsilon^2)(\frac{1}{1-\Delta t})] \left[ e^0 + \left(\frac{T}{1-\Delta t}\right) \left(C_t \Delta t + C_s \dfrac{h^2}{\epsilon^2} \right)  \right]
\end{align}
where the constants $C_p$, $C_t$, and $C_s$ are given by \autoref{eq:constCp} and \autoref{eq:constCtandCs}. The constant $L/\epsilon^2$ is the Lipschitz constant of the peridynamic force $\mathcal{L}^\epsilon(\bu)$ in $L^2$, see \autoref{eq:lipschitzproperty l2}. If the error in initial data is zero then $E^k$ is of the order of $C_t\Delta t + C_s h^2/\epsilon^2$.
\end{theorem}
{\vskip 2mm}
In \autoref{thm:higherregularity}  we have shown that  $\bu,\bv \in C^2([0,T]; W)$ for righthand side $\bb\in C^1([0,T]; W)$.
In \autoref{s:conclusions} we discuss the behavior of the exponential constant appearing in \autoref{thm:convergence}  for evolution times seen in fracture experiments.
Since we are approximating the solution of an ODE on a Banach space the proof of \autoref{thm:convergence} will follow from the Lipschitz continuity of the force $\mathcal{L}^\epsilon(\bu)$ with respect to the $L^2$ norm. The proof is given in the following two sections.

\subsubsection{Truncation error analysis and consistency}
\label{consistent}
The results in this section follows the same steps as in \cite{CMPer-JhaLipton3} and therefore we will just highlight the major steps. We can write the discrete evolution equation for $(\be^k_h(\bu) = \br_h(\bu^k) - \bu^k_h, \be^k_h(\bv) = \br_h(\bv^k) - \bv^k_h)$ as follows
\begin{align}
( \be^{k+1}_h(\bu), \tilde{\bu}) &= (\be^{k}_h(\bu), \tilde{\bu}) + \Delta t (\be^{k+1}_h(\bv), \tilde{\bu} ) + \Delta t (\btau^k_h(\bu),\tilde{\bu}) ,\notag \\
(\be^{k+1}_h(\bv), \tilde{\bu}) &= (\be^k_h(\bv), \tilde{\bu}) + \Delta t (\mathcal{L}^\epsilon(\bu^k_h) - \mathcal{L}^\epsilon(\br_h(\bu^k)), \tilde{\bu})  \notag \\
&\; + \Delta t (\btau^k_h(\bv), \tilde{\bu})+ \Delta t (\bsigma^k_{h}(\bu), \tilde{\bu}), \label{eq:ekevolution}
\end{align}
where consistency error terms $\btau^k_h(\bu), \btau^k_h(\bv), \bsigma^k_h(\bu)$ are given by
\begin{align}
\btau^k_h(\bu) &:= \dfrac{\partial \bu^{k+1}}{\partial t} - \dfrac{\bu^{k+1} - \bu^k}{\Delta t}, \notag \\
\btau^k_h(\bv) &:= \dfrac{\partial \bv^k}{\partial t} - \dfrac{\bv^{k+1} - \bv^k}{\Delta t}, \notag \\
\bsigma^k_{h}(\bu) &:= \mathcal{L}^\epsilon(\br_h(\bu^k)) - \mathcal{L}^\epsilon(\bu^k). \label{eq:consistencyerror}
\end{align}

When $\bu,\bv$ are $C^2$ in time, we easily see that
\begin{align*}
||\btau^k_h(\bu)|| &\leq \Delta t \sup_{t} ||\frac{\partial^2 \bu}{\partial t^2}|| \qquad \text{and} \qquad ||\btau^k_h(\bv)|| \leq \Delta t \sup_{t} ||\frac{\partial^2 \bv}{\partial t^2}||. 
\end{align*}

To estimate $\bsigma^k_{h}(\bu)$, we recall the Lipschitz continuity property of the peridynamic force in the $L^2$ norm, see \autoref{eq:lipschitzproperty l2}. This leads us to 
\begin{align}
||\bsigma^k_{h}(\bu)|| &\leq \dfrac{L}{\epsilon^2} ||\bu^k - \br_h(\bu^k)|| \leq \dfrac{L c}{\epsilon^2} h^2 \sup_{t} ||\bu(t)||_2, \label{eq:consistencyerror peridynamicforce}
\end{align}
where the constant $L$ is defined in \autoref{eq:lipschitzconstant l2}. 

We now state the consistency of this approach.

\begin{lemma}\label{lem:consistency}
\textbf{Consistency}\\
Let $\tau$ be given by
\begin{align}
\tau &:= \sup_{k} \left( ||\btau^k_h(\bu)|| + ||\btau^k_h(\bv)|| + ||\bsigma^k_{h}(\bu)|| \right), \label{eq:tau} 
\end{align}
then the approach is consistent in that
\begin{align}
\tau&\leq C_t \Delta t + C_s \frac{h^2}{\epsilon^2} , \label{eq:totalconsistencyerror}
\end{align}
where
\begin{align}
C_t &:= ||\frac{\partial^2 \bu}{\partial t^2}|| + ||\frac{\partial^2 \bv}{\partial t^2}|| \quad \text{and} \quad C_s := L c \sup_{t} ||\bu(t)||_2.\label{eq:constCtandCs}
\end{align}
\end{lemma}

\subsubsection{Stability analysis} 
\label{stable}
In equation for $\be^k_h(\bu)$, we take $\tilde{\bu} = \be^{k+1}_h(\bu)$. We have
\begin{align*}
||\be^{k+1}_h(\bu)||^2 &= (\be^{k}_h(\bu), \be^{k+1}_h(\bu)) + \Delta t (\be^{k+1}_h(\bv), \be^{k+1}_h(\bu)) + \Delta t (\btau^k_h(\bu), \be^{k+1}_h(\bu)),
\end{align*}
which implies
\begin{align}\label{eq:ekubound}
||\be^{k+1}_h(\bu)|| &\leq ||\be^{k}_h(\bu)|| + \Delta t ||\be^{k+1}_h(\bv)|| + \Delta t ||\btau^{k}_h(\bu)|| .
\end{align}
Similarly, we can show
\begin{align}\label{eq:ekvbound}
||\be^{k+1}_h(\bv)|| &\leq ||\be^{k}_h(\bv)|| + \Delta t ||\mathcal{L}^\epsilon(\bu^k_h) - \mathcal{L}^\epsilon(\br_h(\bu^k))|| \notag \\
& \quad + \Delta t \left( ||\btau^{k}_h(\bv)|| +  ||\bsigma^{k}_{per,h}(\bu)|| \right).
\end{align}

We have from \autoref{eq:lipschitzproperty l2}
\begin{align}\label{eq:stabilityerror peridynamicforce}
||\mathcal{L}^\epsilon(\bu^k_h) - \mathcal{L}^\epsilon(\br_h(\bu^k))|| &\leq \dfrac{L}{\epsilon^2} ||\bu^k_h  - \br_h(\bu^k)|| = \dfrac{L}{\epsilon^2} ||\be^k_h(\bu)||.
\end{align}

After adding \autoref{eq:ekubound} and \autoref{eq:ekvbound}, and substituting \autoref{eq:stabilityerror peridynamicforce}, we get
\begin{align*}
||\be^{k+1}_h(\bu)|| + ||\be^{k+1}_h(\bv)|| &\leq ||\be^k_h(\bu)|| + ||\be^k_h(\bv)|| + \Delta t ||\be^{k+1}_h(\bv)|| + \dfrac{L}{\epsilon^2} \Delta t ||\be^k_h(\bu)|| + \Delta t \tau
\end{align*}
where $\tau$ is defined in \autoref{eq:tau}. Since $e^k = ||\be^k_h(\bu)|| + ||\be^k_h(\bv)||$, we can show, assuming $L/\epsilon^2 \geq 1$,
\begin{align*}
&e^{k+1} \leq e^k + \Delta t e^{k+1} + \Delta t \dfrac{L}{\epsilon^2} e^k + \Delta t \tau \\
\Rightarrow & e^{k+1} \leq \dfrac{1 + \Delta t L/\epsilon^2}{1-\Delta t} e^k + \dfrac{\Delta t}{1 - \Delta t} \tau.
\end{align*}

Substituting for $e^k$ recursively in the equation above, we get
\begin{align*}
e^{k+1} &\leq \left(\dfrac{1 + \Delta t L/\epsilon^2}{1-\Delta t} \right)^{k+1} e^0 + \dfrac{\Delta t}{1 - \Delta t} \tau \sum_{j=0}^k \left(\dfrac{1 + \Delta t L/\epsilon^2}{1-\Delta t} \right)^{k-j}.
\end{align*}
Noting that
\begin{align*}
\dfrac{1 + \Delta t L/\epsilon^2}{1-\Delta t} &= 1 + \frac{(1+L/\epsilon^2)}{1-\Delta t} \Delta t
\end{align*}
and $(1 +a \Delta t )^k \leq \exp[k a\Delta t ] \leq \exp[Ta]$ for $a>0$, we have
\begin{align}
\left(\dfrac{1 + \Delta t L_1/\epsilon^2}{1-\Delta t} \right)^k &\leq \exp[\frac{T(1+L_1/\epsilon^2)}{1-\Delta t}].
\end{align}
This implies
\begin{align*}
e^{k+1} &\leq \exp[\frac{T(1+L/\epsilon^2)}{1-\Delta t}] \left[ e^0 + \dfrac{\Delta t}{1 - \Delta t} \tau \sum_{j=0}^k 1 \right]\\
&\leq \exp[\frac{T(1+L/\epsilon^2)}{1-\Delta t}] \left[ e^0 + \dfrac{k\Delta t}{1 - \Delta t} \tau\right].
\end{align*}
By substituting above equation in \autoref{eq:Ekineq}, we get stability of the scheme:
\begin{lemma}\label{lem:stability}
\textbf{Stability}\\
\begin{align}\label{eq:stab}
E^k &\leq C_p h^2 + \exp[\frac{T(1+L/\epsilon^2)}{1-\Delta t}] \left[ e^0 + \dfrac{k\Delta t}{1 - \Delta t} \tau\right].
\end{align}
\end{lemma}
After taking sup over $k\leq T/\Delta t$ and substituting the bound on $\tau$ from \sref{Lemma}{lem:consistency}, we get the desired result and proof of \autoref{thm:convergence} is complete.

We now consider a stronger notion of stability for the linearized peridynamics model.

\subsection{Linearized peridynamics and energy stability}\label{ss:fullfemstability}
In this section, we linearize the peridynamics model and obtain a CFL like stability condition. For problems where strains are small, the stability condition for the linearized model is expected to apply to the nonlinear model. The slope of peridynamics potential $f$ and $g$ are constant for sufficiently small strain and therefore for small strain the nonlinear model behaves like a linear model. 

In \autoref{eq:nonlocforcetensile}, linearization gives
\begin{align}\label{eq:nonlocforcetensile linear}
\mathcal{L}^{\epsilon}_{T,l}(\bu)(\bx) = \dfrac{2}{\epsilon^{d+1} \omega_d} \int_{H_{\epsilon}(\bx)} \omega(\bx) \omega(\by) J^\epsilon(|\by - \bx|) f''(0) S(\by,\bx;\bu) \be_{\by - \bx} d\by.
\end{align}
The corresponding bilinear form is denoted as $a^\epsilon_{T,l}$ and is given by
\begin{align}\label{eq:operatora Tlinear}
a^\epsilon_{T,l}(\bu, \bv) &= \dfrac{f''(0)}{\epsilon^{d+1} \omega_d} \int_D \int_{D} \omega(\bx) \omega(\by) J^\epsilon(|\by - \bx|) |\by - \bx| S(\by,\bx;\bu) S(\by,\bx;\bv) d\by d\bx. 
\end{align}
Similarly, linearization of $\mathcal{L}^\epsilon_D$ in \autoref{eq:nonlocforcedevia} gives
\begin{align}
\mathcal{L}^\epsilon_{D,l}(\bu)(\bx)=\frac{g''(0)}{\epsilon^{d+2} \omega_d}\int_{H_\epsilon(\bx)} \omega(\bx) \omega(\by) J^\epsilon(|\by-\bx|)\left[\theta(\by,t;\bu)+\theta(\bx,t;\bu)\right]\be_{\by-\bx}\,d\by.
\end{align}
The associated bilinear form is given by
\begin{align}\label{eq:operatora Dlinear}
a^\epsilon_{D,l}(\bu,\bv) = \dfrac{g''(0)}{\epsilon^2} \int_D \omega(\bx) \theta(\bx;\bu) \theta(\by;\bv) d\bx.
\end{align}
The total force after linearization is
\begin{align}
\mathcal{L}^\epsilon_l (\bu)(\bx) = \mathcal{L}^\epsilon_{T,l}(\bu)(\bx) + \mathcal{L}^\epsilon_{D,l}(\bu)(\bx)
\end{align}
and the bilinear operator associated to $\mathcal{L}^\epsilon_l$ is given by
\begin{align}
a^\epsilon_l(\bu,\bv) = a^\epsilon_{T,l}(\bu,\bv) + a^\epsilon_{D,l}(\bu,\bv).
\end{align}
We have
\begin{align*}
(\mathcal{L}^\epsilon_l(\bu), \bv) = - a^\epsilon_l(\bu,\bv).
\end{align*}

We now discuss the stability of the FEM approximation to the linearized problem. Let $\bu^k_{l,h}$ denote the approximate solution satisfying, for $k \geq 1$,
\begin{align}\label{eq:centrallinear}
\left( \dfrac{\bu^{k+1}_{l,h} - 2 \bu^k_{l,h} + \bu^{k-1}_{l,h}}{\Delta t^2}, \tilde{\bu} \right) &= (\mathcal{L}^\epsilon_l(\bu^k_{l,h}), \tilde{\bu} ) + (\bb^k_{h}, \tilde{\bu}), \qquad \forall \tilde{\bu} \in V_h
\end{align}
and, for $k=0$,
\begin{align}\label{eq:centralzerolinear}
\left( \dfrac{\bu^1_{l,h} - \bu^0_{l,h}}{\Delta t^2}, \tilde{\bu}\right) &= \dfrac{1}{2}(\mathcal{L}^\epsilon(\bu^0_{l,h}), \tilde{\bu})+ \dfrac{1}{\Delta t} (\bv^0_{l,h}, \tilde{\bu}) + \dfrac{1}{2}(\bb^0_h, \tilde{\bu}), \qquad \forall \tilde{\bu} \in V_h.
\end{align}

The following notation will be used to define the discrete energy at each time step $k$
\begin{align}
&\overline{\bu}^{k+1}_h := \frac{\bu^{k+1}_h + \bu^k_h}{2}, \: \overline{\bu}^{k}_h := \frac{\bu^k_h + \bu^{k-1}_h}{2}, \notag \\
&\bar{\partial}_t \bu^k_h := \frac{\bu^{k+1}_h - \bu^{k-1}_h}{2 \Delta t}, \: \bar{\partial}_t^+ \bu^k_h :=  \frac{\bu^{k+1}_h - \bu^{k}_h}{\Delta t}, \: \bar{\partial}_t^- \bu^k_h :=  \frac{\bu^{k}_h - \bu^{k-1}_h}{\Delta t}.
\end{align}
We also define
\begin{align*}
\bar{\partial}_{tt} \bu^k_h &:= \dfrac{\bu^{k+1}_h - 2\bu^k_h + \bu^{k-1}_h}{\Delta t^2} = \dfrac{\bar{\partial}^+_t \bu^k_h - \bar{\partial}^-_t \bu^k_h}{\Delta t}.
\end{align*}

We introduce the discrete energy associated with $\bu^k_{l,h}$ at time step $k$ as follows
\begin{align*}
\mathcal{E}(\bu^k_{l,h}) &:= \frac{1}{2} \left[  ||\bar{\partial}^+_t \bu^k_{l,h}||^2 - \frac{\Delta t^2}{4} a^\epsilon_l(\bar{\partial}^+_t \bu^k_{l,h}, \bar{\partial}^+_t \bu^k_{l,h}) + a^\epsilon_l(\overline{\bu}^{k+1}_{l,h} , \overline{\bu}^{k+1}_{l,h}) \right]
\end{align*}

Following [Theorem 4.1, \cite{CMPer-Karaa}], the stability of central difference scheme is given by
\begin{theorem}\label{thm:cflcondition}
\textbf{Energy Stability of the Central difference approximation of linearized peridynamics}\\
Let $\bu^k_{l,h}$ be the approximate solution of \autoref{eq:centrallinear} and \autoref{eq:centralzerolinear}. In the absence of body force $\bb(t) = 0$ for all $t$, if $\Delta t$ satisfies the CFL like condition
\begin{align}
\frac{\Delta t^2}{4} \sup_{\bu\in V_h \setminus \{\bzero\}} \dfrac{a^\epsilon_l(\bu,\bu)}{(\bu,\bu)} \leq 1,
\end{align}
then the discrete energy is positive and we have the stability
\begin{align}\label{eq:wcoercive}
\mathcal{E}(\bu^k_{l,h})= \mathcal{E}(\bu^{0}_{l,h}) .
\end{align}
\end{theorem}

We skip the proof of above theorem as it is straight-forward extension of Theorem 5.2 in \cite{CMPer-JhaLipton3}.

\section{Numerical experimants}\label{s:numerical}
In this section, we present numerical simulations that are consistent with the theoretical a-priori bound on the convergence rate. We also compare the peridynamic energy of the material softening zone and the  classic Griffith's fracture energy of linear elastic fracture mechanics.

We consider Plexiglass at room temperature and specify the density $\rho = 1200 \,kg/m^3$, bulk modulus $K = 25 \, GPa$, Poisson's ratio $\nu = 0.245$, and critical energy release rate $G_c = 500 \,Jm^{-2}$. The pairwise interaction and the hydrostatic interaction are characterized by potentials $f(r) = c (1-\exp [-\beta r^2])$ and $g(r) = \bar{C} r^2/2$ respectively. Here we have used a quadratic hydrostatic interaction potential. The influence function is $J(r) = 1-r$. Since pairwise potential $f$ is symmetric for positive and negative strain the critical strain is given by $S_c(\by, bx) = \frac{\pm\bar{r}}{\sqrt{\by - \bx}}$, where $\pm\bar{r}$ is the inflection point of $f(r)$ given by $\bar{r} = \frac{1}{\sqrt{\beta}}$. 
Following equations 94, 95, and 97 of  \cite{CMPer-Lipton4}, the relation between peridynamic material parameters and Lam\'e constants $(\lambda, \mu)$ and critical energy release rate $G_c$ can be written as (for 2-d)
\begin{align}\label{eq:materialrelation}
c = \frac{\pi G_c}{4M_J}, \qquad \beta = \frac{4\mu}{C M_J}, \qquad \bar{C} = \frac{2(\lambda - \mu)}{M_J^2},
\end{align}
where $M_J$ is given by
\begin{align*}
M_J &= \int_0^1 J(r) r^2 dr = \frac{1}{12}.
\end{align*}
By solving \autoref{eq:materialrelation}, we get $c = 4712.4$, $\bar{C}=-1.7349\times 10^{11}$,  $\beta = 1.5647\times 10^{8}$. 


\begin{figure}
\centering
\includegraphics[scale=0.15]{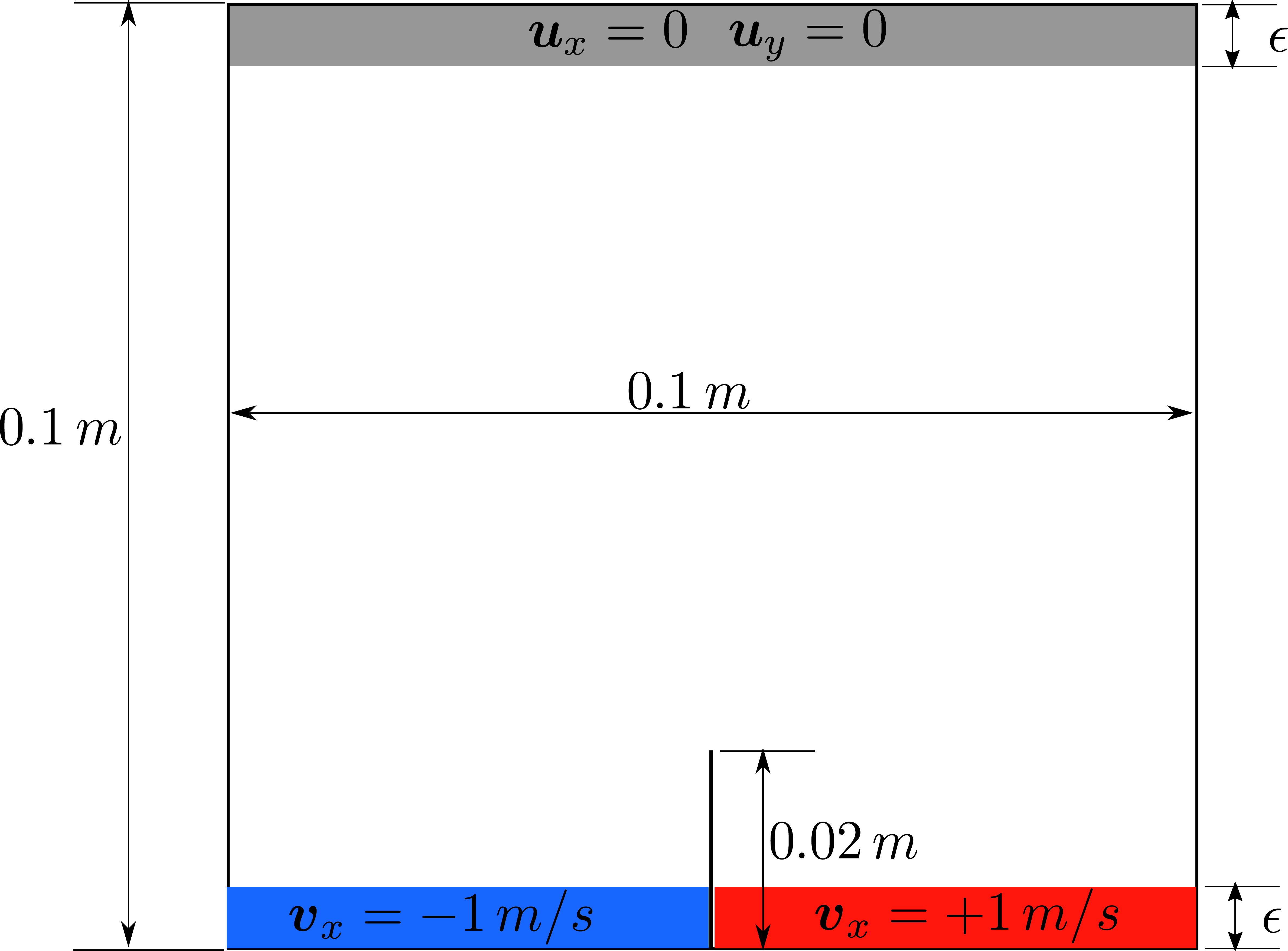}
\caption{Material domain $D = [0,0.1\, m]^2$ with crack of length $0.02\,m$. The x-component and y-component of displacement are fixed along a collar of thickness equal to the horizon on top. On the bottom the velocity $\bv_x = \pm 1\, m/s$ along x-direction is specified on either side of the crack to make the crack propagate upwards.}\label{fig:setup crack prop}
\end{figure}

We consider a 2-d domain $D=[0,0.1\, m]^2$ (with unit thickness in third direction) with vertical crack of length $0.02\, m$. The boundary conditions are described in \autoref{fig:setup crack prop}.  The simulation time is $T = 40\, \mu s$ and the time step is $\Delta t = 0.004\,\mu s$. We consider two horizons $8\, mm$ and $4\, mm$. We run simulations for mesh sizes $h = 2, 1, 0.5\, mm $. We consider the central difference time discretization described by \autoref{eq:central} on a uniform mesh consisting of linear triangle elements. Second order quadrature approximation is used in the simulation for each triangle element. To reduce the load on memory and to avoid matrix-vector multiplication at each time step, we approximate the mass matrix by diagonal mass matrix using lumping (row-sum) technique. Suppose exact mass matrix is $\bolds{M} = [m_{ij}]$ where $m_{ij}$ is the element of $\bolds{M}$ corresponding to $i^{\text{th}}$ row and $j^{\text{th}}$ column, then we approximate $\bolds{M}$ by diagonal matrix $\hat{\bolds{M}} = [\hat{m}_{ij}]$ where $\hat{m}_{ii} = \sum_{j} m_{ij}$ and $\hat{m}_{ij} = 0$ if $ j\neq i$.

\subsection{Convergence rate}
To compute convergence rate numerically we proceed as follows: consider a fixed horizon $\epsilon$ and three different mesh sizes $h_1, h_2, h_3$ such that $r = h_1/h_2 = h_2/h_3$. Let $\bu_1,\bu_2,\bu_3$ be approximate solutions corresponding to meshes of size  $h_1,h_2,h_3$, and let $\bu$ be the exact solution. We write the error as $||\bu_h - \bu|| =C h^\alpha$ for some constant $C$ and $\alpha>0$, to get
\begin{align*}
\log( ||\bu_1 - \bu_2||) &= C + \alpha \log h_2, \\
\log( ||\bu_2 - \bu_3||) &= C + \alpha \log h_3.
\end{align*}
From above two equations, it is easy to see that the rate of convergence $\alpha$ is 
\begin{align}\label{eq:rate formula}
\dfrac{\log( ||\bu_1 - \bu_2||) - \log( ||\bu_2 - \bu_3||)}{\log(r)}.
\end{align}

%

\begin{figure}[h]
\centering
\includegraphics[scale=0.4]{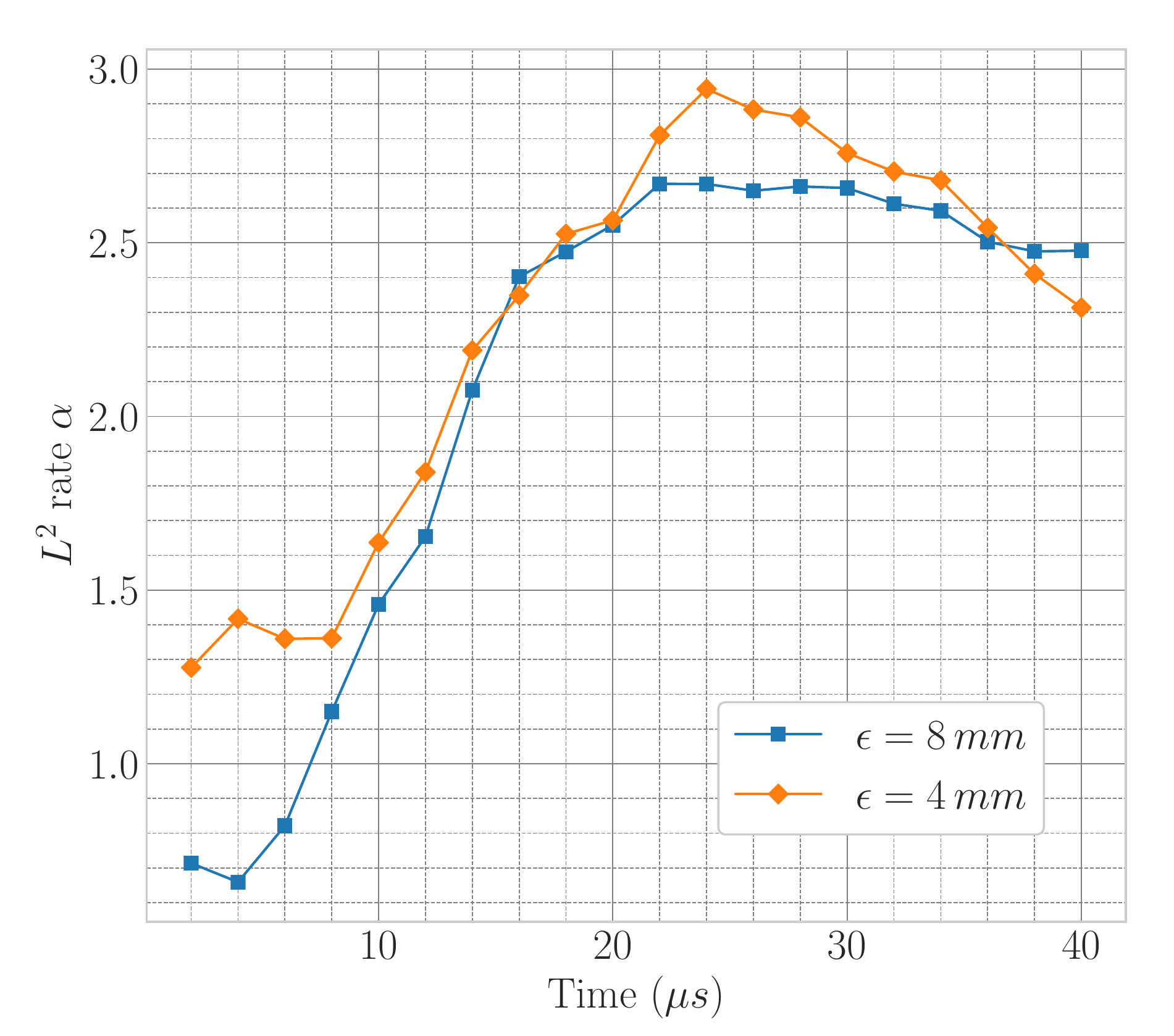}
\caption{Convergence rate at different times for two horizons. For both horizons $\epsilon = 4, 8\, mm$, the three meshes of size $h = 2, 1, 0.5\, mm$ were considered to compute the convergence rate.}\label{fig:rate}
\end{figure}

The convergence result for horizons $\epsilon = 8\,mm$ and $\epsilon = 4\, mm$ is shown in \autoref{fig:rate}. In the simulation we have considered second order approximation of integration using quadrature points. The simulations show a rate of convergence that agrees with the a-priori estimates given in \autoref{thm:convergence}.

\subsection{Fracture energy of crack zone}\label{sec: softzone}
The extent of damage at material point $\bx$ is given by the function $Z(\bx)$
\begin{align}\label{eq:damage}
Z(\bx) &= \max_{\by \in H_\epsilon(\bx) \cap D} \frac{S(\by,\bx;\bu)}{S^+_c}.
\end{align}
The crack zone is defined as set of material points which have $Z > 1$. We compute the peridynamic energy of crack zone and compare it with the Griffith's fracture energy. For a crack of length $l$, the Griffith's fracture energy (G.E.) will be $G.E. = G_c \times l$. The peridynamic fracture energy (P.E.) associated with the material softening zone is given by 
\begin{align*}
P.E. = &\int_{\substack{\bx \in D,\\ Z(\bx) \geq 1}} \left[ \frac{1}{\epsilon^d \omega_d} \int_{H_\epsilon(\bx)} |\by - \bx| \mathcal{W}^\epsilon(S(\by,\bx;\bu))\,d\by \right] d\bx\\
&+\int_{\substack{\bx \in D,\\ Z(\bx) \geq 1}}  \mathcal{V}^\epsilon(\theta(\bx,t;\bu))\,d\bx
\end{align*}
where $\mathcal{W}^\epsilon(S(\by,\bx;\bu))$ is the bond-based potential, see \autoref{eq:tensilepot} and $\mathcal{V}^\epsilon(\theta(\bx,t;\bu))$ is the hydrostatic interaction potential, see \autoref{eq:hydropot}. 

\begin{figure}[h]
\centering
\includegraphics[scale=0.35]{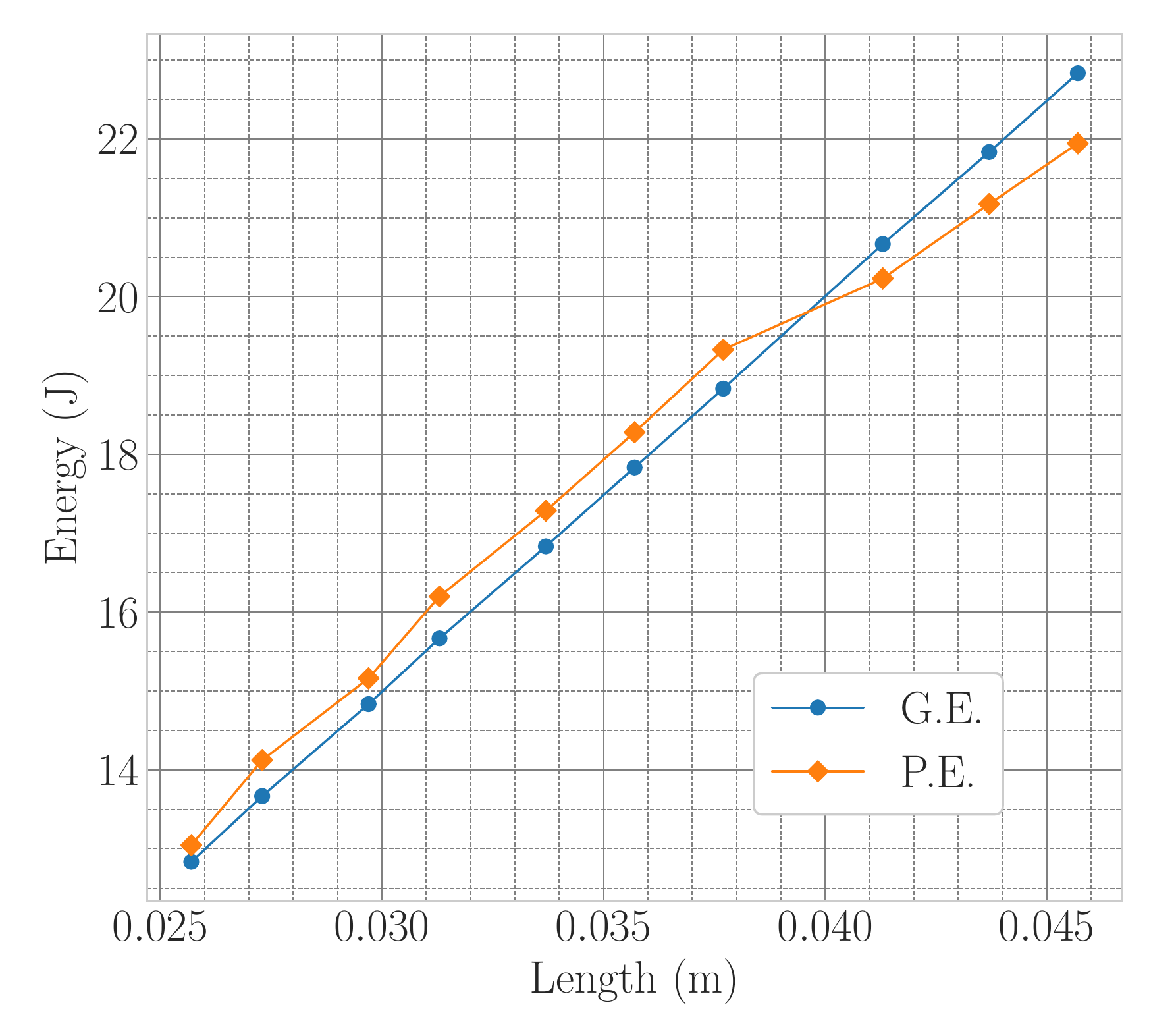}
\caption{Peridynamic energy and Griffith's energy as a function of crack length.}\label{fig:crackenergyvslength}
\end{figure}
\begin{figure}[h]
\centering
\includegraphics[scale=0.3]{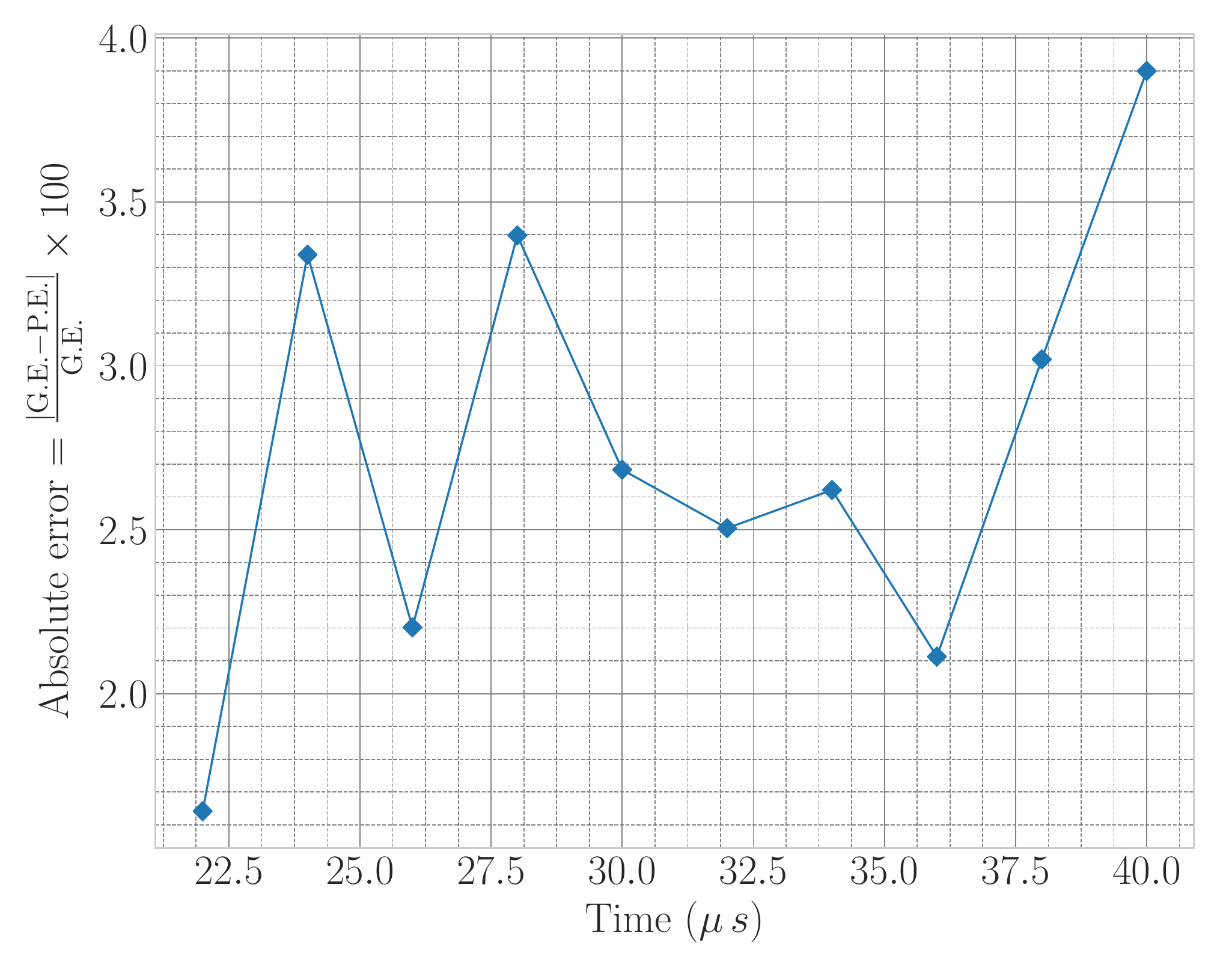}
\caption{Error between Peridynamic energy and Griffith's energy at different times.}\label{fig:crackenergyerror}
\end{figure}
\begin{figure}
    \centering
    \begin{subfigure}{.47\linewidth}
        \centering
        \scalebox{1.0}[1.0]{\includegraphics[scale=0.18]{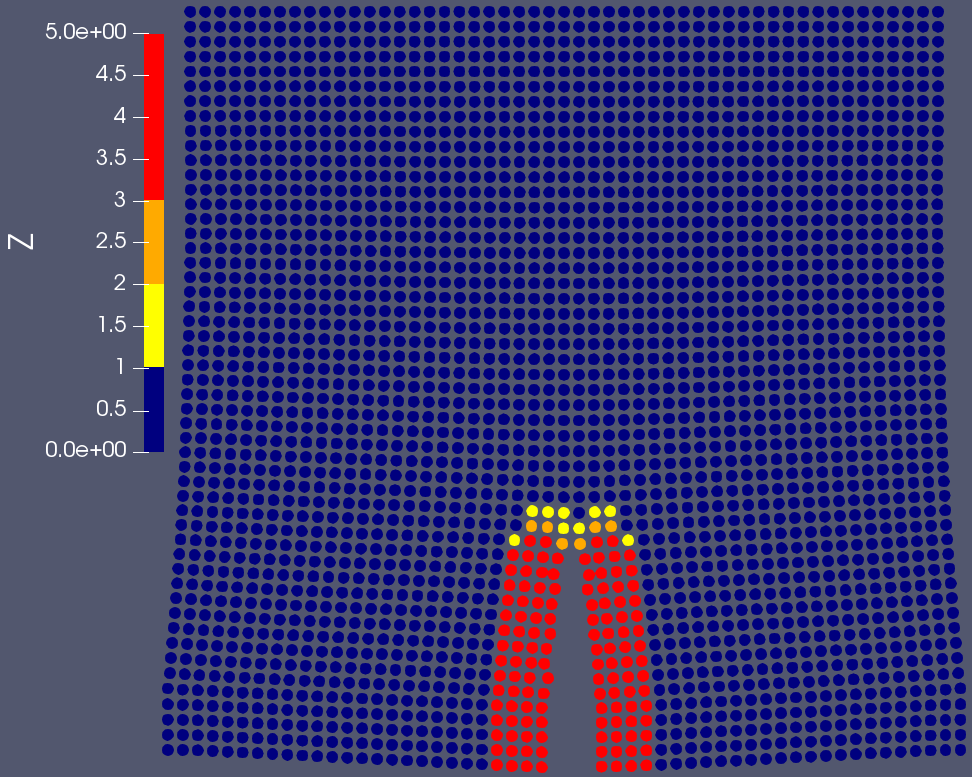}}
        \caption{$t = 30\, \mu s$}
    \end{subfigure}
    \begin{subfigure}{.47\linewidth}
        \centering
        \scalebox{1.0}[1.0]{\includegraphics[scale=0.18]{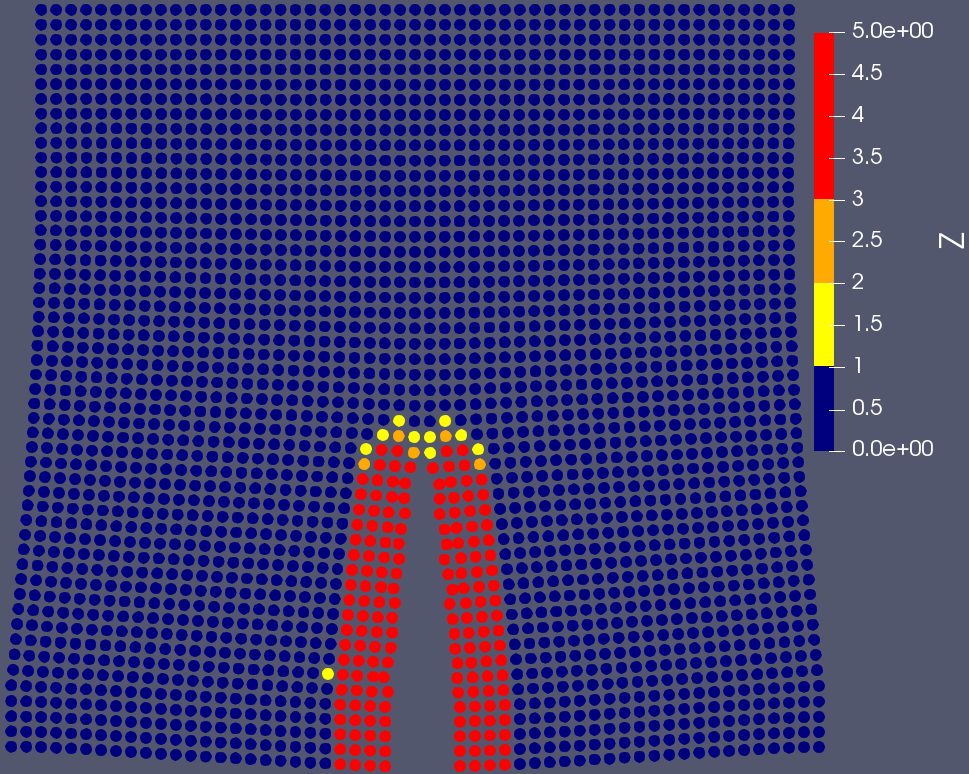}}
        \caption{$t = 40\, \mu s$}
    \end{subfigure}
    
   \caption{Color plot of damage function $Z$ on deformed material domain at time $t=30\,\mu s$ and $40\,\mu s$. Dark blue represents undamaged material $Z<1$, $Z\approx 1$ is yellow at crack tip, red is softening material. Here, the displacements are scaled by $100$ and damage function is cut off at $5$ to highlight the crack zone.}
    \label{fig:damageplot}
\end{figure}
\begin{figure}
    \centering
    \begin{subfigure}{.47\linewidth}
        \centering
        \scalebox{1.7}[1.7]{\includegraphics[scale=0.135]{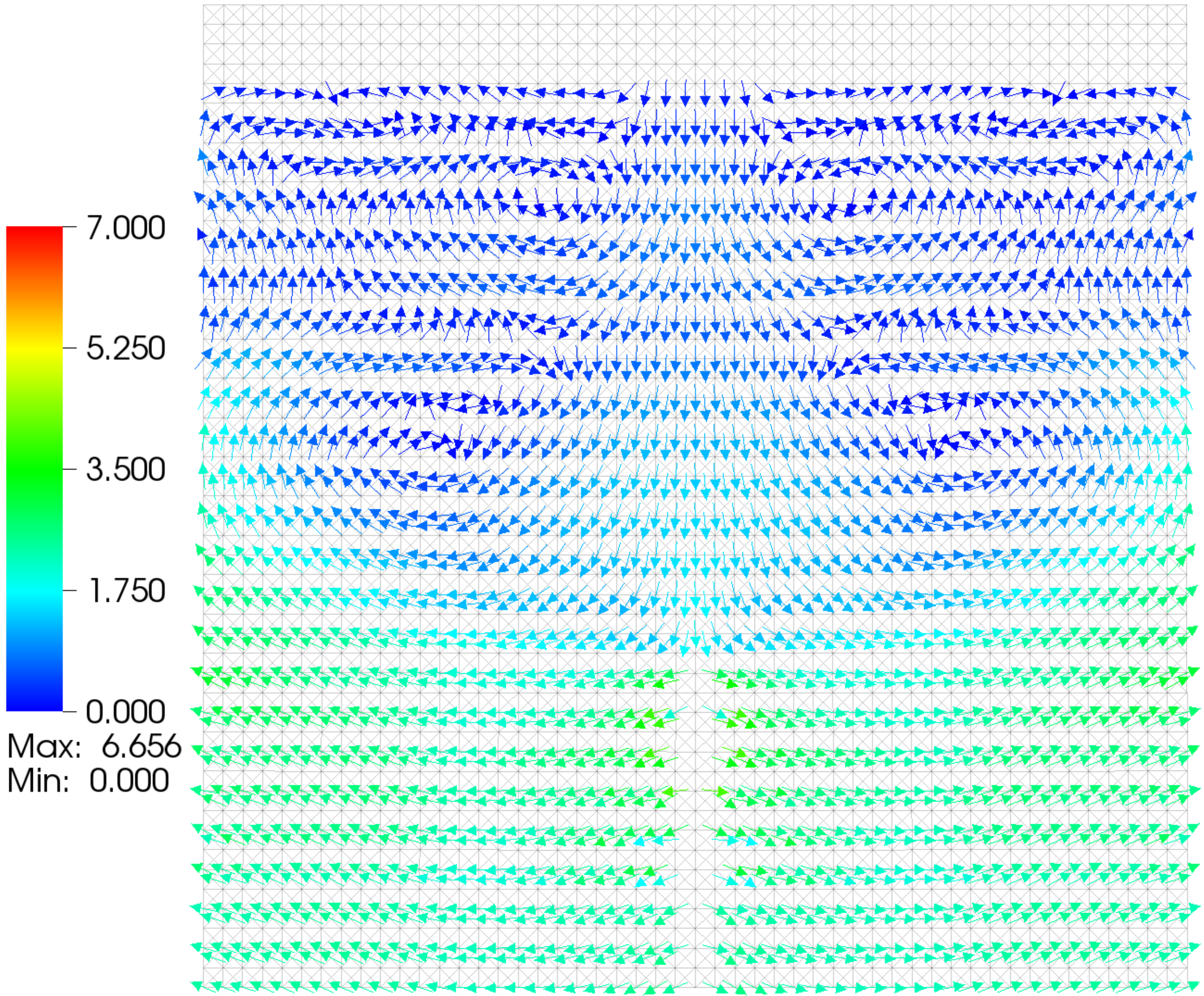}}
        \caption{$t = 30\, \mu s$}
    \end{subfigure}
    \begin{subfigure}{.47\linewidth}
        \centering
        \scalebox{1.7}[1.7]{\includegraphics[scale=0.135]{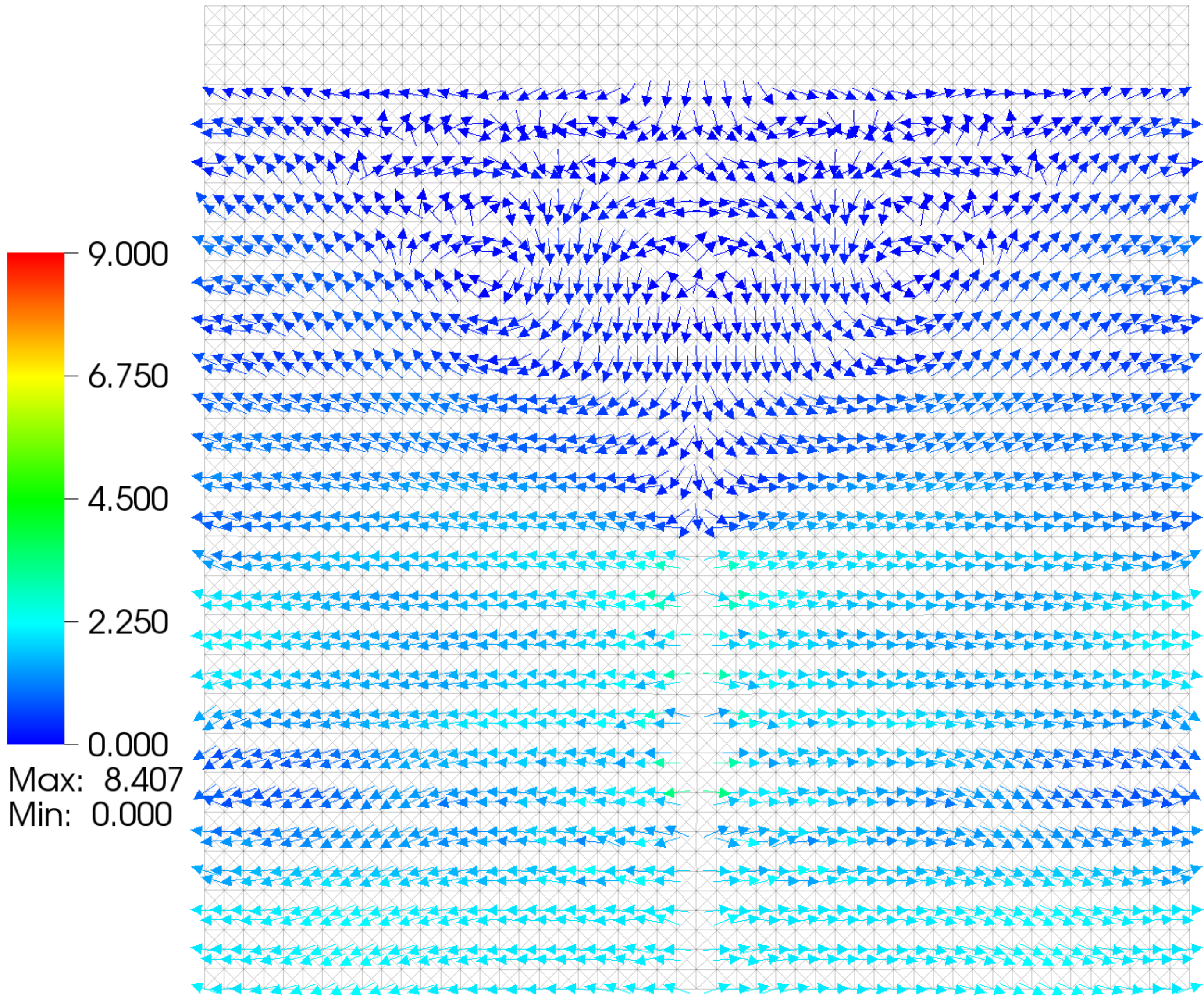}}
        \caption{$t = 40\, \mu s$}
    \end{subfigure}
    
   \caption{Velocity profile.}
    \label{fig:velplot}
\end{figure}
\begin{figure}
    \centering
    \begin{subfigure}{.47\linewidth}
        \centering
        \scalebox{1.0}[1.0]{\includegraphics[scale=0.14]{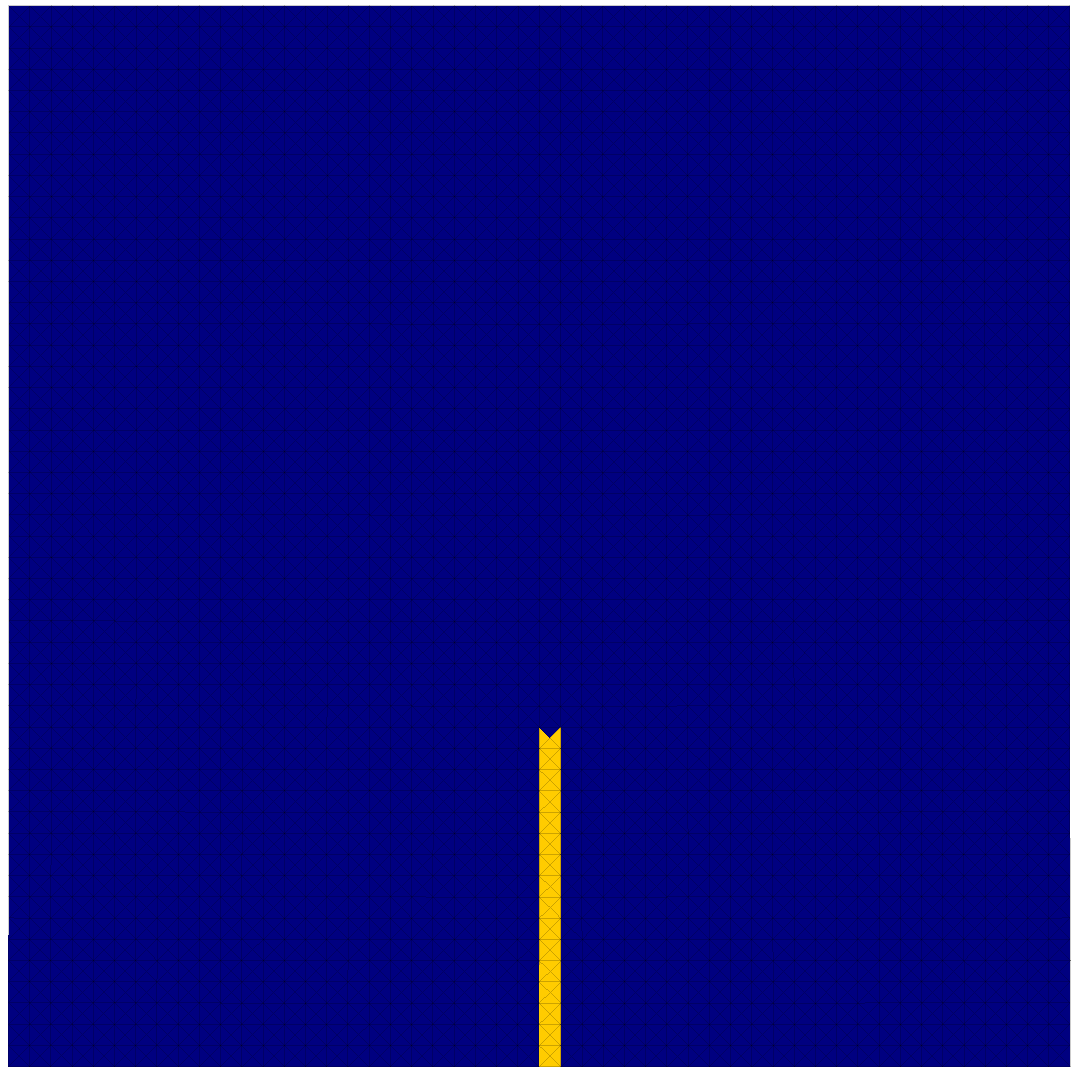}}
        \caption{$t = 30\, \mu s$}
    \end{subfigure}
    \begin{subfigure}{.47\linewidth}
        \centering
        \scalebox{1.0}[1.0]{\includegraphics[scale=0.14]{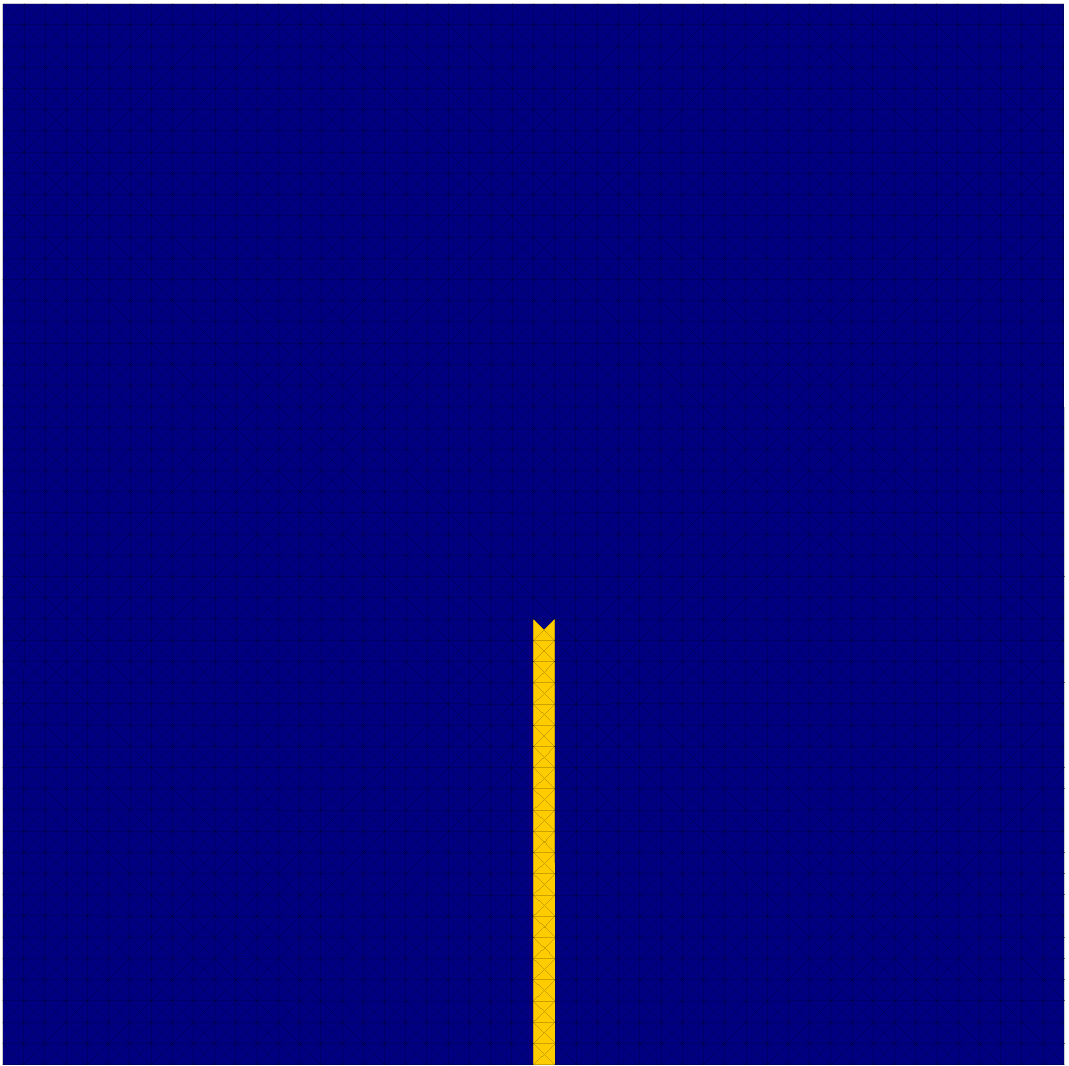}}
        \caption{$t = 40\, \mu s$}
    \end{subfigure}
    
   \caption{Magnitude of the $xx$ component of strain $\nabla \bu + \nabla \bu ^T$. The region for which the magnitude of the strain is greater than a multiple of the critical strain is the yellow region}
    \label{fig:uxxplot}
\end{figure}

In \autoref{fig:crackenergyvslength} classical fracture energy and peridynamic fracture energy is shown at different crack length. The error in both energies at different times is shown in \autoref{fig:crackenergyerror}. The agreement between two energies is good. The damage profile at time $30\,\mu s$ and $40\,\mu s$ is shown in \autoref{fig:damageplot}. At each node, the damage function $Z$ is computed by treating edges between mesh nodes as bonds. 
{In addition to the damage plots, we show velocity profile at $30\, \mu s$ and $40\, \mu s$ in \autoref{fig:velplot}.} In \autoref{fig:uxxplot} we show the plot of the $xx$ component of symmetric gradient of the displacement. Here the region for which the magnitude of the strain is greater than a multiple of the critical strain is the yellow region. It is seen that the high strain region surrounds  the crack.

As the crack is propagating vertically it is seen that the high strain region is next to the crack i

\section{Lipschitz continuity of peridynamic force and higher temporal regularity of solutions}\label{s:proofs}
In this section, we prove \autoref{thm:lipschitzproperty} and \autoref{thm:higherregularity}. Here $\bu \in W\subset H^2(D;\bbR^d)$ and the $||\bu||_2$ norm is given by
\begin{align}
||\bu||_2 = ||\bu|| + ||\nabla \bu|| + ||\nabla^2 \bu||.
\end{align}

\subsection{Proof of Lipschitz continuity with respect to the $\Vert \cdot\Vert_2$ norm}
We assume that the potential function $f$ satisfies $C^f_i < \infty $ for $i=0,1,2,3,4$. Recall that $C^f_0 = \sup_r |f(r)|$ and $C^f_i = \sup_r |f^{(r)}(r)|$ for $i=1,...,4$. $C^g_i$ is defined similarly for $i=0,1,...,4$. If the potential function $g$ is a convex-concave function then we can assume $C^g_i < \infty$ for $i=0,1,2,3,4$. In what follows, we will prove the \autoref{thm:lipschitzproperty} for convex-concave type $g$. If $g$ is a purely a quadratic function the proof follows easily using only a subset of the estimates proved in this section. 

Let $\bu,\bv \in W$. Using the triangle inequality, we get
\begin{align}
||\mathcal{L}^\epsilon(\bu) - \mathcal{L}^\epsilon(\bv)||_2 &\leq ||\mathcal{L}^\epsilon_T(\bu) - \mathcal{L}^\epsilon_T(\bv)||_2 + ||\mathcal{L}^\epsilon_D(\bu) - \mathcal{L}^\epsilon_D(\bv)||_2,
\end{align}
where $\mathcal{L}^\epsilon_T$ and $\mathcal{L}^\epsilon_D$ is given by \autoref{eq:nonlocforcetensile} and \autoref{eq:nonlocforcedevia}. 

We first write the peridynamic force $\mathcal{L}^\epsilon_T(\bu)(\bx)$ as follows
\begin{align}
&\mathcal{L}^\epsilon_T(\bu)(\bx) \notag \\
&=\frac{2}{\epsilon^{d+1} \omega_d}\int_{H_\epsilon(\bx)} \omega(\bx) \omega(\by) \frac{J^\epsilon(|\by-\bx|)}{\sqrt{|\by-\bx|}} f'(\sqrt{|\by-\bx|}S(\by,\bx;\bu))\be_{\by-\bx}\,d\by,
\end{align}
where we substituted $\partial_S f(\sqrt{|\by - \bx|} S(\by,\bx;\bu)) = \sqrt{|\by - \bx|}f'(\sqrt{|\by - \bx|} S(\by,\bx;\bu))$. The form of the peridynamic force described above is the same as the one given in [Section 6, \cite{CMPer-JhaLipton3}]. 
We apply Theorem 3.1 in \cite{CMPer-JhaLipton3}] to show
\begin{align}
||\mathcal{L}^\epsilon_T(\bu) - \mathcal{L}^\epsilon_T(\bv)||_2 &\leq \dfrac{L_1(1+ (||\bu||_2 + ||\bv||_2) + (||\bu||_2 + ||\bv||_2)^2)}{\epsilon^3} ||\bu - \bv||_2 \notag \\
&\leq \dfrac{L_1(1+ ||\bu||_2 + ||\bv||_2)^2}{\epsilon^3} ||\bu - \bv||_2
\end{align}
and
\begin{align}
||\mathcal{L}^\epsilon_T(\bu)||_2 &\leq \dfrac{L_2(||\bu||_2 + ||\bu||_2^2)}{\epsilon^{5/2}}.
\end{align}

Next we analyze $||\mathcal{L}^\epsilon_D(\bu) - \mathcal{L}^\epsilon_D(\bv)||_2$. We define new terms to simplify the calculations. For $\bxi \in H_1(\bzero)$, we  set
\begin{align}\label{eq:notationsproof}
s_\bxi &= \epsilon |\bxi|, \: \be_\bxi = \frac{\bxi}{|\bxi|}, \notag \\
\omega_\bxi(\bx) &= \omega(\bx + \epsilon \bxi) \omega(\bx), \notag \\
\bubar_\bxi(\bx) &= \bu(\bx + \epsilon \bxi) - \bu(\bx), \notag \\
(\bu - \bv)(\bx) &= \bu(\bx) - \bv(\bx).
\end{align}
Similar notations hold if we exchange $\bx$, $\bxi \in H_1(\bzero)$, and $\bu\in W$  by $\by$, $\bolds{\eta} \in H_1(\bzero)$, and  $\bv\in W$ respectively We will also encounter various moments of the influence function $J$ therefore we define following moments
\begin{align}\label{eq:barJ}
\bar{J}_\alpha = \frac{1}{\omega_d} \int_{H_1(\bzero)} J(|\bxi|) |\bxi|^{-\alpha} d\bxi, \qquad \text{ for } \alpha \in \bbR. 
\end{align}
Recall that $J(|\bxi|) = 0$ for $\bxi \notin H_1(\bzero)$ and $0\leq J(|\bxi|) \leq M$ for $\bxi\in H_1(\bzero)$. The boundary function $\omega$ is assumed to satisfy
\begin{align}
\sup_{\bx} |\nabla \omega(\bx)| < \infty ,\qquad \sup_{\bx} |\nabla^2 \omega(\bx)| < \infty.
\end{align}
We choose finite constants $C_{\omega_1}$ and $C_{\omega_2}$ such that
\begin{align}
|\nabla \omega_\bxi(\bx)| &\leq C_{\omega_1}, \quad |\nabla \omega(\bx)| \leq C_{\omega_1}, \notag \\
|\nabla^2 \omega_\bxi(\bx)| &\leq C_{\omega_2}, \quad |\nabla^2 \omega(\bx)| \leq C_{\omega_2}.
\end{align}

We now collect the following estimates which will be used to estimate $||\mathcal{L}^\epsilon_D(\bu) - \mathcal{L}^\epsilon_D(\bv)||_2$. 

{\vskip 2mm}
\begin{lemma}\label{lem:estimateontheta}
Let $\bu,\bv \in W$, for any $\betanew \in H_1(\bzero)$ and $\delta \leq 2\epsilon$, we have
\begin{align}
\sup_{\bx \in D} |\theta(\bx;\bu)| &\leq 2 C_{e_1} \bar{J}_0 ||\bu||_2, \label{eq:est1} \\
\int_D |\theta(\bx + \delta \betanew; \bu)|^2 d\bx &\leq 4 \bar{J}^2_0 ||\bu||_2^2, \label{eq:est2} \\
\int_D |\nabla \theta(\bx + \delta \betanew; \bu)|^2 d\bx &\leq 8 \bar{J}^2_0 (1+C_{\omega_1})^2 ||\bu||_2^2, \label{eq:est3} \\
\int_D |\theta(\bx + \delta \betanew; \bu - \bv)|^2 \, |\nabla \theta(\bx + \delta \betanew;\bv)|^2 d\bx &\leq 32 \bar{J}_0^4 (1+C_{\omega_1})^2 ||\bv||_2^2 \, ||\bu - \bv||_2^2, \label{eq:est4} \\
\int_D |\nabla \theta(\bx + \delta \betanew; \bu)|^4 d\bx &\leq 128 \bar{J}^4_0 (C_{e_2}^2 + C_{e_1} C^2_{\omega_1})^2 ||\bu||_2^4, \label{eq:est5} \\
\int_D |\theta(\bx + \delta \betanew; \bu - \bv)|^2 \,  |\nabla \theta(\bx + \delta \betanew; \bv)|^4 d\bx &\leq 512 \bar{J}_0^6 C_{e_1}^2 (C_{e_2}^2 + C_{e_1}C^2_{\omega_1})^2 ||\bu - \bv||_2^2 \, ||\bv||_2^4, \label{eq:est6} \\
\int_D |\nabla^2 \theta(\bx + \delta \betanew; \bu)|^2 d\bx &\leq 16 \bar{J}^2_0 (1+2C_{\omega_1} + C_{\omega_2})^2 ||\bu||_2^2, \label{eq:est7} .
\end{align}
Here $\nabla$ in all the equations above is with respect to $\bx$. The constants $C_{e_1}, C_{e_2}$ are the constants associated to the Sobolev embedding property of space $H^2(D;\bbR^d)$, see \autoref{eq:sob embedd 1} and \autoref{eq:sob embedd 2}.
\end{lemma}
{\vskip 2mm}
\begin{proof}
Using the notation given in \autoref{eq:notationsproof}, we write $\theta(\bx;\bu)$ as
\begin{align}
\theta(\bx;\bu) &= \frac{1}{\omega_d} \int_{H_1(\bzero)} \omega(\bx+\epsilon \bxi) J(|\bxi|) \bar{\bu}_\bxi (\bx) \cdot \be_\bxi d\bxi.
\end{align}
On noting that  $|\bubar_\bxi(\bx)| \leq 2||\bu||_{\infty}$ and $||\bu||_\infty \leq C_{e_1} ||\bu||_2$, we easily see that
\begin{align}
|\theta(\bx;\bu)| &\leq \bar{J}_0 2 ||\bu||_{\infty} \leq 2 C_{e_1} \bar{J}_0 ||\bu||_2.
\end{align}
In the rest of the proof, we will let $\by = \bx + \delta \betanew$, where $0\leq \delta \leq 2\epsilon$ and $\betanew \in H_1(\bzero)$.

To show \autoref{eq:est2}, we first introduce an important identity which will be used frequently. Let $p(\bxi)$ be some function of $\bxi$, and $\alpha, C \in \bbR$  then
\begin{align}\label{eq:symminequality}
&\bigg\vert \frac{C}{\omega_d} \int_{H_1(\bzero)} \frac{J(|\bxi|)}{|\bxi|^\alpha} p(\bxi) d\bxi \bigg\vert^2 \notag \\
&= \left(\frac{C}{\omega_d} \right)^2 \int_{H_1(\bzero)} \int_{H_1(\bzero)} \frac{J(|\bxi|)}{|\bxi|^\alpha} \frac{J(|\boldsymbol{\eta}|)}{|\boldsymbol{\eta}|^\alpha} p(\bxi) p(\boldsymbol{\eta}) d\bxi d\boldsymbol{\eta}  \notag \\
&\leq \left(\frac{C}{\omega_d} \right)^2 \int_{H_1(\bzero)} \int_{H_1(\bzero)} \frac{J(|\bxi|)}{|\bxi|^\alpha} \frac{J(|\boldsymbol{\eta}|)}{|\boldsymbol{\eta}|^\alpha} \frac{p(\bxi)^2 + p(\boldsymbol{\eta})^2}{2} d\bxi d\boldsymbol{\eta} \notag \\
&=  C^2 \frac{\bar{J}_\alpha}{\omega_d } \int_{H_1(\bzero)} \frac{J(|\bxi|)}{|\bxi|^\alpha} p(\bxi)^2 d\bxi,
\end{align}
where we used the inequality $ab \leq a^2/2 + b^2/2$ in the first step, and definition of $\bar{J}_\alpha$ and symmetry of terms in second step.

{From the expression of $\theta(\by;\bu)$ we can show
\begin{align}
\int_D |\theta(\by;\bu)|^2 d\bx &\leq 4 \bar{J}^2_0 ||\bu||^2 \leq 4 \bar{J}^2_0 ||\bu||^2_2.
\end{align}
}

We now prove the bound \autoref{eq:est3}. Taking the gradient of $\theta(\by;\bu)$, with respect to $\bx$, noting that $\by = \bx + \delta \betanew$, we get
\begin{align}\label{eq:gradtheta}
\nabla \theta (\by;\bu) &= \frac{1}{\omega_d} \int_{H_1(\bzero)} J(|\bxi|) \omega(\by + \epsilon \bxi) (\nabla \bubar_\bxi(\by))^T \be_\bxi d\bxi \notag \\
&\quad + \frac{1}{\omega_d} \int_{H_1(\bzero)} J(|\bxi|) \nabla \omega(\by + \epsilon \bxi) \bubar_\bxi(\by) \cdot \be_\bxi d\bxi.
\end{align}
We can show using inequality \autoref{eq:symminequality} and the estimates $\int_D |\nabla \bubar_\bxi(\by)|^2 d\bx \leq 4 ||\nabla \bu||^2 \leq 4 ||\bu||_2^2$, $| \nabla \omega(\by + \epsilon \bxi)| \leq C_{\omega_1}$, $\int_D |\bubar_\bxi(\by)|^2 d\bx \leq 4 ||\bu||_2^2$, to conclude
\begin{align}
\int_D |\nabla \theta(\by ; \bu)|^2 d\bx &\leq \dfrac{2\bar{J}_0}{\omega_d} \int_{H_1(\bzero)} J(|\bxi|) 4 ||\bu||_2^2 d\bxi + \dfrac{2\bar{J}_0}{\omega_d} \int_{H_1(\bzero)} J(|\bxi|) 4 C_{\omega_1}^2||\bu||_2^2 d\bxi \notag \\
&= 8 \bar{J}^2_0 (1+C_{\omega_1}^2) ||\bu||_2^2 \leq 8 \bar{J}^2_0 (1+C_{\omega_1})^2 ||\bu||_2^2.
\end{align}

We now show \autoref{eq:est4}. We will use \autoref{eq:est1} and \autoref{eq:est3}, and proceed as follows
\begin{align}
& \int_D |\theta(\by; \bu - \bv)|^2 \, |\nabla \theta(\by;\bv)|^2 d\bx \notag  \\
&\leq \left(\sup_{\by} |\theta(\by; \bu - \bv)| \right)^2 \int_D |\nabla \theta(\by;\bv)|^2 d\bx \notag  \\
&\leq 32 \bar{J}_0^4 (1+C_{\omega_1})^2 ||\bv||^2_2 ||\bu - \bv||^2_2.
\end{align}

{To prove \autoref{eq:est5} we note expression of $\nabla \theta(\by;\bu)$ in \autoref{eq:gradtheta} and inequality $(a+b)^4 \leq 8 a^4 + 8 b^4$ and \autoref{eq:symminequality} to get
\begin{align}
|\nabla \theta(\by ;\bu)|^4 &\leq \dfrac{64 \bar{J}^3_0}{\omega_d} \int_{H_1(\bzero)} J(|\bxi|) (|\nabla \bu(\by+\epsilon \bxi)|^4 + |\nabla \bu(\by)|^4)  d\bxi \notag \\
&\quad + \dfrac{64 C^4_{\omega_1} \bar{J}^3_0}{\omega_d} \int_{H_1(\bzero)} J(|\bxi|) (|\bu(\by+\epsilon \bxi)|^4 + |\bu(\by)|^4) d\bxi.
\end{align}}
Application of  Fubini's theorem gives
\begin{align}
&\int_D |\nabla \theta(\by ;\bu)|^4 d\bx \notag \\
&\leq \dfrac{64 \bar{J}^3_0}{\omega_d} \int_{H_1(\bzero)} J(|\bxi|) \left( \int_D (|\nabla \bu(\by+\epsilon \bxi)|^4 + |\nabla \bu(\by)|^4) d\bx \right)  d\bxi \notag \\
&\quad + \dfrac{64 C^4_{\omega_1} \bar{J}^3_0}{\omega_d} \int_{H_1(\bzero)} J(|\bxi|) \left( \int_D (|\bu(\by+\epsilon \bxi)|^4 + |\bu(\by)|^4) d\bx \right) d\bxi \notag \\
&\leq \dfrac{64 \bar{J}^3_0}{\omega_d} \int_{H_1(\bzero)} J(|\bxi|) \left( 2 ||\nabla \bu||^4_{L^4(D;\bbR^{d\times d})} \right)  d\bxi \notag \\
&\quad + \dfrac{64 C^4_{\omega_1} \bar{J}^3_0}{\omega_d} \int_{H_1(\bzero)} J(|\bxi|) \left( ||\bu||^2_\infty \int_D (|\bu(\by+\epsilon \bxi)|^2 + |\bu(\by)|^2) d\bx \right) d\bxi \notag \\
&\leq 128 \bar{J}^4_0  ||\nabla \bu||^4_{L^4(D;\bbR^{d\times d})}  + \dfrac{64 C^4_{\omega_1} \bar{J}^3_0}{\omega_d} \int_{H_1(\bzero)} J(|\bxi|) \left( ||\bu||^2_\infty 2 ||\bu||^2 \right) d\bxi \notag \\
&\leq 128 \bar{J}^4_0  ||\nabla \bu||^4_{L^4(D;\bbR^{d\times d})}  + 128 C^4_{\omega_1} \bar{J}^4_0 ||\bu||^2_\infty  \,||\bu||^2.
\end{align}
Using the Sobolev embedding property, $||\bu||_\infty \leq C_{e_1}||\bu||_2$ and $||\nabla \bu||_{L^4} \leq C_{e_2} ||\bu||_2$, we obtain
\begin{align}
\int_D |\nabla \theta(\by ;\bu)|^4 d\bx &\leq 128 \bar{J}^4_0 (C_{e_2}^4 + C_{e_1}^2 C^4_{\omega_1}) ||\bu||^4_2 \leq 128 \bar{J}^4_0 (C_{e_2}^2 + C_{e_1} C^2_{\omega_1})^2 ||\bu||^4_2 \,.
\end{align}

The estimate \autoref{eq:est6} follows by combining estimates \autoref{eq:est1} and \autoref{eq:est5}. It now remains to show \autoref{eq:est7}. From expression of $\nabla \theta(\by;\bu)$ in \autoref{eq:gradtheta}, we have
\begin{align}\label{eq:gradgradtheta}
\nabla^2 \theta(\by;\bu) &= \frac{1}{\omega_d} \int_{H_1(\bzero)} J(|\bxi|) \omega(\by + \epsilon \bxi) \nabla^2 (\bubar_\bxi(\by) \cdot \be_\bxi) d\bxi \notag \\
&\quad + \frac{1}{\omega_d} \int_{H_1(\bzero)} J(|\bxi|) ((\nabla (\bubar_\bxi(\by))^T \be_\bxi) \dyad \nabla \omega(\by + \epsilon \bxi) d\bxi \notag \\
&\quad + \frac{1}{\omega_d} \int_{H_1(\bzero)} J(|\bxi|) \nabla \omega(\by + \epsilon \bxi) \dyad ((\nabla (\bubar_\bxi(\by))^T \be_\bxi) d\bxi \notag \\
&\quad + \frac{1}{\omega_d} \int_{H_1(\bzero)} J(|\bxi|) \nabla^2 \omega(\by + \epsilon \bxi) \bubar_\bxi(\by) \cdot \be_\bxi d\bxi.
\end{align}
{Using the equation above we can show
\begin{align}
\int_D |\nabla^2 \theta(\by;\bu)|^2 d\bx &\leq \frac{3 \bar{J}_0}{\omega_d} \int_{H_1(\bzero)} J(|\bxi|) \left( \int_D |\nabla^2 \bubar_\bxi(\by)|^2 d\bx \right) d\bxi \notag \\
&\quad + \frac{12 C^2_{\omega_1} \bar{J}_0}{\omega_d} \int_{H_1(\bzero)} J(|\bxi|) \left( \int_D |\nabla \bubar_\bxi(\by)|^2 d\bx \right) d\bxi \notag \\
&\quad + \frac{3C^2_{\omega_2} \bar{J}_0}{\omega_d} \int_{H_1(\bzero)} J(|\bxi|) \left( \int_D | \bubar_\bxi(\by)|^2 d\bx \right) d\bxi.
\end{align}}
The terms $|\bubar_\bxi(\by)|^2$, $|\nabla \bubar_\bxi(\by)|^2$, and $|\nabla^2 \bubar_\bxi(\by)|^2$ are bounded by $2 (|\bu(\by + \epsilon\bxi)|^2 + |\bu(\by)|^2)$, $2 (|\nabla \bu(\by + \epsilon\bxi)|^2 + |\nabla \bu(\by)|^2)$, and $2 (|\nabla^2 \bu(\by + \epsilon\bxi)|^2 + |\nabla^2 \bu(\by)|^2)$ respectively. Therefore, we have
\begin{align}
\int_D |\nabla^2 \theta(\by;\bu)|^2 d\bx &\leq \left[ 3 \bar{J}^2_0 + 12 C^2_{\omega_1} \bar{J}^2_0 + 3 C^2_{\omega_2} \bar{J}^2_0\right] 4 ||\bu||^2_2 \notag \\
&\leq 16 \bar{J}^2_0 (1+C_{\omega_2} + 2C_{\omega_1})^2 ||\bu||^2_2,
\end{align}
and this completes the proof of lemma. 
\end{proof}
{\vskip 2mm}

{\vskip 2mm}
\textbf{Estimating $||\mathcal{L}^\epsilon_D(\bu) - \mathcal{L}^\epsilon_D(\bv)||$: }We apply the  notation described in \autoref{eq:notationsproof}, and write $\mathcal{L}^\epsilon_D(\bu)(\bx)$ as follows
\begin{align}\label{eq:forcedevia}
\mathcal{L}^\epsilon_D (\bu)(\bx) &= \frac{1}{\epsilon^2 \omega_D} \int_{H_1(\bzero)} \omega_\bxi(\bx) J(|\bxi|) [g'(\theta(\bx + \epsilon\bxi;\bu)) + g'(\theta(\bx;\bu))] \be_\bxi d\bxi.
\end{align}
{Using the formula above and from the expression for $\theta$ we can easily show
\begin{align}\label{LD}
||\mathcal{L}^\epsilon_D (\bu) - \mathcal{L}^\epsilon_D (\bv)|| &\leq \dfrac{L_1}{\epsilon^2} ||\bu - \bv||_2,
\end{align}
where $L_1 = 4 C^g_2 \bar{J}_0^2$.
}

{\vskip 2mm}
\textbf{Estimating $||\nabla \mathcal{L}^\epsilon_D(\bu) - \nabla \mathcal{L}^\epsilon_D(\bv)||$: }Taking the gradient of \autoref{eq:forcedevia} gives
\begin{align}\label{eq:gradforcedevia}
\nabla \mathcal{L}^\epsilon_D (\bu) (\bx) &= \dfrac{1}{\epsilon^2 \omega_d} \int_{H_1(\bzero)} J(|\bxi|) \omega_\bxi(\bx) \be_\bxi \dyad [\nabla g'(\theta(\bx + \epsilon \bxi; \bu)) + \nabla g'(\theta(\bx; \bu))] d\bxi \notag \\
&\quad + \dfrac{1}{\epsilon^2 \omega_d} \int_{H_1(\bzero)} J(|\bxi|) \be_\bxi \dyad \nabla \omega_\bxi(\bx) [g'(\theta(\bx + \epsilon \bxi; \bu)) + g'(\theta(\bx; \bu))] d\bxi \notag \\
&=: G_1(\bu)(\bx) + G_2(\bu)(\bx),
\end{align}
where we have denoted the first and second term as $G_1(\bu)(\bx)$ and $G_2(\bu)(\bx)$ for convenience. On using the triangle inequality we get
\begin{align*}
||\nabla \mathcal{L}^\epsilon_D (\bu) - \nabla \mathcal{L}^\epsilon_D (\bv)|| &\leq ||G_1(\bu) - G_1(\bv)|| + ||G_2(\bu) - G_2(\bu)||. 
\end{align*}

From the expression of $G_1(\bu)$, we have
\begin{align*}
|G_1(\bu)(\bx) - G_1(\bv)(\bx) | &\leq \dfrac{1}{\epsilon^2 \omega_d} \int_{H_1(\bzero)} J(|\bxi|) (|\nabla g'(\theta(\bx + \epsilon \bxi; \bu)) - \nabla g'(\theta(\bx + \epsilon \bxi; \bv))| \notag \\
&\qquad \qquad + |\nabla g'(\theta(\bx; \bu)) - \nabla g'(\theta(\bx; \bv))|) d\bxi.
\end{align*}
Let 
\begin{align}\label{eq:p1}
p_1(\by) &:= |\nabla g'(\theta(\by; \bu)) - \nabla g'(\theta(\by; \bv))|
\end{align}
{and we get
\begin{align}\label{eq:estinterm2}
||G_1(\bu) - G_1(\bv)||^2 &\leq \left(\dfrac{1}{\epsilon^2}\right)^2 \dfrac{2\bar{J}_0}{ \omega_d} \int_{H_1(\bzero)} J(|\bxi|) \left( \int_D (p_1(\bx + \epsilon \bxi)^2 + p_1(\bx)^2) d\bx \right) d\bxi.
\end{align}}
Note that
\begin{align*}
\nabla g'(\theta(\bx + \epsilon \bxi; \bu)) &= g''(\theta(\bx + \epsilon\bxi;\bu)) \nabla \theta(\bx + \epsilon\bxi;\bu).
\end{align*}
Therefore, from \autoref{eq:p1},
\begin{align*}
p_1(\by) &= |g''(\theta(\by;\bu)) \nabla \theta(\by;\bu) - g''(\theta(\by;\bv)) \nabla \theta(\by;\bv)| \notag \\
&\leq C^g_2 |\nabla \theta(\by;\bu) - \nabla \theta(\by;\bv)| +  C^g_3 |\theta(\by;\bu) - \theta(\by;\bv)| \, |\nabla \theta(\by;\bv)| \notag \\
&= C^g_2 |\nabla \theta(\by;\bu - \bv)| +  C^g_3 |\theta(\by;\bu-\bv)| \, |\nabla \theta(\by;\bv)|,
\end{align*}
where we have added and subtracted $g''(\theta(\by;\bu)) \nabla \theta(\by;\bv)$ and used the fact that $g''(r) \leq C^g_2$ and $|g''(r_1) - g''(r_2)| \leq C^g_3 |r_1 - r_2|$. We use the estimate on $p_1$ and proceed as follows
\begin{align*}
\int_D p_1(\by)^2 d\bx &\leq 2 (C^g_2)^2 \int_D |\nabla \theta(\by;\bu - \bv)|^2 d\bx \notag \\
&\quad + 2(C^g_3)^2 \int_D  |\theta(\by;\bu-\bv)|^2 \, |\nabla \theta(\by;\bv)|^2 d\bx,
\end{align*}
where we denote $\bx+\epsilon \bxi$ as $\by$. We apply inequality \autoref{eq:est3} and \autoref{eq:est4} of \sref{Lemma}{lem:estimateontheta} to obtain
\begin{align*}
\int_D p_1(\by)^2 d\bx &\leq 2 (C^g_2)^2 8 \bar{J}^2_0 (1+C_{\omega_1})^2 ||\bu - \bv||_2^2\notag \\
&\quad  + 2(C^g_3)^2 32 \bar{J}_0^4 (1+C_{\omega_1})^2 ||\bv||_2^2 \, ||\bu - \bv||_2^2 \notag \\
&\leq L_2 (1+ ||\bv||_2)^2 ||\bu - \bv||_2^2, 
\end{align*}
where we have grouped all the constant factors together and denote their product by $L_2$. Substituting these estimates in to \autoref{eq:estinterm2} gives
\begin{align}\label{eq:estinterm3}
||G_1(\bu) - G_1(\bv)||^2 &\leq \dfrac{4L_2 \bar{J}_0^2 }{\epsilon^4} (1+ ||\bv||_2)^2 ||\bu - \bv||_2^2 \notag \\
\Rightarrow ||G_1(\bu) - G_1(\bv)|| &\leq \dfrac{L_3(1+ ||\bv||_2)}{\epsilon^2} ||\bu - \bv||_2,
\end{align}
where we have introduced the new constant $L_3$. 

{The formula for $G_2(\bu)$ is similar to $\mathcal{L}^\epsilon_D(\bu)$ and therefore we have 
\begin{align*}
||G_2(\bu) - G_2(\bv)|| &\leq \dfrac{C_{\omega_1} L_1}{\epsilon^2} ||\bu - \bv||_2. 
\end{align*}}
Collecting results, we have shown
\begin{align}\label{gradLD}
||\nabla \mathcal{L}^\epsilon_D (\bu) - \nabla \mathcal{L}^\epsilon_D (\bv)|| &\leq \dfrac{L_4(1+||\bv||_2)}{\epsilon^2} ||\bu - \bv||_2,
\end{align}
where we have introduced new constant $L_4$.

{\vskip 2mm}
\textbf{Estimating $||\nabla^2 \mathcal{L}^\epsilon_D(\bu) - \nabla^2 \mathcal{L}^\epsilon_D(\bv)||$: }Taking the gradient of \autoref{eq:gradforcedevia}, gives
\begin{align}\label{eq:gradgradforcedevia}
\nabla^2 \mathcal{L}^\epsilon_D(\bu)(\bx) &= \dfrac{1}{\epsilon^2 \omega_d} \int_{H_1(\bzero)} \omega_\bxi(\bx) J(|\bxi|) \be_\bxi \dyad [\nabla^2 g'(\theta(\bx+\epsilon\bxi;\bu)) + \nabla^2 g'(\theta(\bx;\bu))] d\bxi \notag \\
&\quad + \dfrac{1}{\epsilon^2 \omega_d} \int_{H_1(\bzero)} J(|\bxi|) \be_\bxi \dyad [\nabla g'(\theta(\bx+\epsilon\bxi;\bu)) + \nabla g'(\theta(\bx;\bu))] \dyad \nabla \omega_\bxi(\bx) d\bxi \notag \\
&\quad + \dfrac{1}{\epsilon^2 \omega_d} \int_{H_1(\bzero)} J(|\bxi|) \be_\bxi \dyad \nabla \omega_\bxi(\bx) \dyad [\nabla g'(\theta(\bx+\epsilon\bxi;\bu)) + \nabla g'(\theta(\bx;\bu))] d\bxi \notag \\
&\quad + \dfrac{1}{\epsilon^2 \omega_d} \int_{H_1(\bzero)} J(|\bxi|) \be_\bxi \dyad \nabla^2 \omega_\bxi(\bx) [g'(\theta(\bx+\epsilon\bxi;\bu)) + g'(\theta(\bx;\bu))] d\bxi \notag \\
&=: H_1(\bu)(\bx) + H_2(\bu)(\bx) + H_3(\bu)(\bx) + H_4(\bu)(\bx).
\end{align}
It is easy to see that estimate on $||H_2(\bu) - H_2(\bv)||$ and $||H_3(\bu) - H_3(\bv)||$ is similar to the estimate for $||G_1(\bu) - G_1(\bv)||$. Thus, from \autoref{eq:estinterm3}, we have
\begin{align}\label{eq:estinterm4}
||H_2(\bu) - H_2(\bv)|| + ||H_3(\bu) - H_3(\bv)|| &\leq \dfrac{2C_{\omega_1} L_3 (1+||\bv||_2) }{\epsilon^2} ||\bu - \bv||_2. 
\end{align}
Also the estimate for $||H_4(\bu) - H_4(\bv)||$ is similar to the estimate for $||G_2(\bu) - G_2(\bv)||$ and we conclude
\begin{align}\label{eq:estinterm5}
||H_4(\bu) - H_4(\bv)|| &\leq \dfrac{C_{\omega_2} L_1 }{\epsilon^2} ||\bu - \bv||_2.
\end{align}

We now work on $||H_1(\bu) - H_1(\bv)||$. From expression of $H_1(\bu)(\bx)$ in \autoref{eq:gradgradforcedevia}, we can easily get following
\begin{align}
|H_1(\bu)(\bx) - H_1(\bv)(\bx)| &\leq \dfrac{1}{\epsilon^2 \omega_d} \int_{H_1(\bzero)} J(|\bxi|) (|\nabla^2 g'(\theta(\bx + \epsilon \bxi; \bu)) - \nabla^2 g'(\theta(\bx + \epsilon \bxi; \bv))| \notag \\
&\qquad \qquad + |\nabla^2 g'(\theta(\bx; \bu)) - \nabla^2 g'(\theta(\bx; \bv))|) d\bxi.
\end{align}
Let $p_2(\by)$, where $\by = \bx + \epsilon \bxi$ and $\nabla$ is with respect to $\bx$, is given by 
\begin{align}
p_2(\by) &:= |\nabla^2 g'(\theta(\by; \bu)) - \nabla^2 g'(\theta(\by; \bv))|.
\end{align}
{We then have
\begin{align}\label{eq:estinterm6}
||H_1(\bu) - H_1(\bv)||^2 &\leq \left(\dfrac{1}{\epsilon^2}\right)^2 \dfrac{2\bar{J}_0}{ \omega_d} \int_{H_1(\bzero)} J(|\bxi|) \left( \int_D (p_2(\bx + \epsilon \bxi)^2 + p_2(\bx)^2) d\bx \right) d\bxi.
\end{align}
}
Note that 
\begin{align*}
\nabla^2 g'(\theta(\by; \bu)) = g'''(\theta(\by;\bu)) \nabla \theta(\by;\bu) \dyad \nabla \theta(\by;\bu) + g''(\theta(\by;\bu)) \nabla^2 \theta(\by;\bu).
\end{align*}
We add and subtract terms to the equation above to get
\begin{align*}
&\nabla^2 g'(\theta(\by; \bu) - \nabla^2 g'(\theta(\by; \bv) \notag \\
&= g'''(\theta(\by;\bu)) [ \nabla \theta(\by;\bu) \dyad \nabla \theta(\by;\bu) - \nabla \theta(\by;\bv) \dyad \nabla \theta(\by;\bv)] \notag \\
&\quad + [g'''(\theta(\by;\bu)) - g'''(\theta(\by;\bv))] \nabla \theta(\by;\bv) \dyad \nabla \theta(\by;\bv) \notag \\
&\quad + g''(\theta(\by;\bu)) [\nabla^2 \theta(\by;\bu) -  \nabla^2 \theta(\by;\bv)] \notag \\
&\quad + [g''(\theta(\by;\bu)) - g''(\theta(\by;\bv))] \nabla^2 \theta(\by;\bv).
\end{align*}
Using inequalities $|g''(r)| \leq C^g_2$, $|g'''(r)|\leq C^g_3$, $|g''(r_1) - g''(r_2)| \leq C^g_3 |r_1 - r_2|$, $|g'''(r_1) - g'''(r_2)| \leq C^g_4 |r_1 - r_2|$, and $|\ba\dyad \ba - \bc\dyad \bc| \leq (|\ba| + |\bc|) |\ba - \bc|$, and the fact that $\theta(\by;\bu) - \theta(\by;\bv) = \theta(\by;\bu - \bv)$, we have
\begin{align*}
p_2(\by) &\leq C^g_3 |\nabla \theta(\by;\bu)| \, |\nabla \theta(\by;\bu - \bv)| + C^g_3 |\nabla \theta(\by;\bv)|\, |\nabla \theta(\by;\bu - \bv)| \notag \\
&\quad + C^g_4 | \theta(\by;\bu - \bv)| \, |\nabla \theta(\by;\bv)|^2 \notag \\
&\quad + C^g_2 |\nabla^2 \theta(\by;\bu-\bv)| + C^g_3 |\theta(\by;\bu - \bv)|\, |\nabla^2 \theta(\by;\bv)|.
\end{align*}
Taking the square of the above equation and using $(\sum_{i=1}^5 a_i)^2 \leq 5 \sum_{i=1}^5 a^2_i$ gives
\begin{align*}
\int_D p_2(\by)^2 d\bx &\leq 5(C^g_3)^2 \int_D |\nabla \theta(\by;\bu)|^2 \, |\nabla \theta(\by;\bu - \bv)|^2 d\bx \notag \\
&\quad + 5(C^g_3)^2 \int_D |\nabla \theta(\by;\bv)|^2 \, |\nabla \theta(\by;\bu - \bv)|^2 d\bx \notag \\
&\quad + 5 (C^g_4)^2 \int_D | \theta(\by;\bu - \bv)|^2 \, |\nabla \theta(\by;\bv)|^4 d\bx \notag \\
&\quad + 5 (C^g_2)^2 \int_D |\nabla^2 \theta(\by;\bu-\bv)|^2 d\bx \notag \\
&\quad + 5(C^g_3)^2 \int_D |\theta(\by;\bu - \bv)|^2 \, |\nabla^2 \theta(\by;\bv)|^2 d\bx \notag \\
&=: I_1 + I_2 + I_3 + I_4 + I_5.
\end{align*}

We now estimate each term using \sref{Lemma}{lem:estimateontheta} as follows. 
Applying H\"older inequality and inequality \autoref{eq:est5} of \sref{Lemma}{lem:estimateontheta} to get
\begin{align*}
I_1 &\leq 5(C^g_3)^2 \left(\int_D |\nabla \theta(\by;\bu)|^4 d\bx  \right)^{1/2} \left(\int_D |\nabla \theta(\by;\bu - \bv)|^4 d\bx  \right)^{1/2}  \notag \\
&\leq 640 (C^g_3)^2 \bar{J}_0^4 (C^2_{e_2} + C_{e_1}C^2_{\omega_1})^2 ||\bu||_2^2\, ||\bu - \bv||_2^2.
\end{align*}
Similarly,
\begin{align*}
I_2 &\leq 640 (C^g_3)^2 \bar{J}_0^4 (C^2_{e_2} + C_{e_1}C^2_{\omega_1})^2 ||\bv||_2^2\, ||\bu - \bv||_2^2.
\end{align*}
Using \autoref{eq:est6} of \sref{Lemma}{lem:estimateontheta}, we get
\begin{align*}
I_3 &\leq 2560 (C^g_4)^2 \bar{J}_0^6 C_{e_1}^2 (C^2_{e_2} + C_{e_1}C^2_{\omega_1})^2 ||\bv||_2^4\, ||\bu - \bv||_2^2.
\end{align*}
For $I_4$, we use inequality \autoref{eq:est7} to get
\begin{align*}
I_4 &\leq 80 (C^g_2)^2 \bar{J}_0^2 (1+2C_{\omega_1} + C_{\omega_2})^2 ||\bu - \bv||_2^2. 
\end{align*}
In $I_5$, we use \autoref{eq:est1} and \autoref{eq:est7} to get
\begin{align*}
I_5 &\leq 320 \bar{J}^4_0 C_{e_1}^2 (1+2C_{\omega_1} + C_{\omega_2})^2 ||\bv||_2^2\, ||\bu - \bv||_2^2.
\end{align*}
After collecting results, we can find a constant $L_5$ such that we have
\begin{align}\label{eq:estinterm7}
\int_D p_2(\by)^2 d\bx &\leq L_5^2 \,(1+ (||\bu||_2 + ||\bv||_2) + (||\bu||_2 + ||\bv||_2)^2)^2 \,||\bu - \bv||_2^2.
\end{align}

We  substitute \autoref{eq:estinterm7} into \autoref{eq:estinterm6} to show
\begin{align}
||H_1(\bu) - H_1(\bv)||^2 &\leq \frac{4 L^2_5 \bar{J}_0^2 (1+ (||\bu||_2 + ||\bv||_2) + (||\bu||_2 + ||\bv||_2)^2)^2}{\epsilon^4} ||\bu - \bv||_2^2 \notag \\
\Rightarrow ||H_1(\bu) - H_1(\bv)|| &\leq \dfrac{L_6(1+ (||\bu||_2 + ||\bv||_2) + (||\bu||_2 + ||\bv||_2)^2) }{\epsilon^2} ||\bu - \bv||_2 \notag \\
&\leq \dfrac{L_6(1+ ||\bu||_2 + ||\bv||_2)^2 }{\epsilon^2} ||\bu - \bv||_2 
\end{align}
where we have introduced the new constant $L_6$. 

We combine the estimates on $H_1,H_2,H_3,H_4$, introducing a new constant $L_7$, and get
\begin{align}\label{eq:estinterm8}
||\nabla^2 \mathcal{L}^\epsilon_D(\bu) - \nabla^2 \mathcal{L}^\epsilon_D(\bv)||
 &\leq \dfrac{L_7 (1+ ||\bu||_2 + ||\bv||_2)^2 }{\epsilon^2} ||\bu - \bv||_2. 
\end{align}
On adding the estimates  \autoref{LD}, \autoref{gradLD}, \autoref{eq:estinterm8} it is evident that the proof of \autoref{thm:lipschitzproperty} is complete. 

\subsection{Proof of higher temporal regularity}
In this section we prove that the peridynamic evolutions have higher regularity in time for body forces that that are differentiable in time.
To see this we take the time derivative of \autoref{eq:equationofmotion} to get a second order differential equation in time for $\bv=\dot{\bu}$ given by
\begin{align}\label{eq:equationofmotionv}
\rho \partial^2_{tt} \bv(\bx,t) &= Q(\bv(t); \bu(t))(\bx) + \dot{\bb}(\bx,t),
\end{align}
where $Q(\bv;\bu)$ is an operator that depends on the solution $\bu$ of  \autoref{eq:equationofmotion} and acts on $\bv$. It is given by, 
\begin{align}\label{eq:Q}
Q(\bv;\bu)(\bx) &= Q_T(\bv;\bu)(\bx) + Q_D(\bv;\bu)(\bx),\hbox{    $\forall \bx\in D$,}
\end{align}
where 
\begin{align}\label{eq:QT}
Q_T(\bv;\bu)(\bx) =\frac{2}{\epsilon^d \omega_d}\int_{H_\epsilon(\bx)} &\omega(\bx) \omega(\by) \frac{J^\epsilon(|\by-\bx|)}{\epsilon|\by-\bx|} \notag \\
&\partial^2_{SS} f(\sqrt{|\by-\bx|}S(\by,\bx,t;\bu)) S(\by,\bx,t;\bv) \be_{\by-\bx}\,d\by,
\end{align}
and
\begin{align}\label{eq:QD}
Q_D(\bv;\bu)(\bx) =\frac{1}{\epsilon^d \omega_d} &\int_{H_\epsilon(\bx)} \omega(\bx) \omega(\by) \frac{J^\epsilon(|\by-\bx|)}{\epsilon^2} \notag \\
&\;\left[\partial^2_{\theta \theta} g(\theta(\by,t;\bu)) \theta(\by,t;\bv) +\partial^2_{\theta\theta} g(\theta(\bx,t;\bu)) \theta(\bx,t;\bv) \right]\be_{\by-\bx}\,d\by.
\end{align}
Clearly, for $\bu$ fixed, the form $Q(\bv;\bu)$ acts linearly on $\bv$ which implies that the equation for $\bv$ is a linear nonlocal equation. The linearity of $Q(\bv;\bu)$ implies Lipschitz continuity for $\bv \in W$ as stated below.

\begin{theorem}\label{thm:lipschitzpropertyQ}
\textbf{Lipschitz continuity of $Q$}\\
Let $\bu \in W$ be any given field. Then for all $\bv,\bw \in W$, we have
\begin{align}
||Q(\bv;\bu) - Q(\bw;\bu)||_2 &\leq \dfrac{L_8(1 + ||\bu||_2)^2}{\epsilon^3} ||\bv-\bw||_2\label{eq:lipschitzpropertyQ}
\end{align}
where the constant $L_8$ does not depend on $\bu,\bv,\bw$. This gives for all $\bv\in W$ the upper bound,
\begin{align}
||Q(\bv;\bu)||_2 &\leq \dfrac{L_8(1 + ||\bu||_2)^2}{\epsilon^3} ||\bv||_2 \label{eq:upperboundQ}.
\end{align}
\end{theorem}
The proof follows the same steps used in the proof of \autoref{thm:lipschitzproperty}. 

If $\bu$ is a peridynamic solution such that $\bu \in C^2(I_0;W)$ then we have for all $t \in I_0$ the inequality
\begin{align}\label{eq:upperboundQtime}
||Q(\bv;\bu(t))||_2 &\leq \dfrac{L_8(1 + \sup_{s \in I_0} ||\bu(s)||_2)^2}{\epsilon^3} ||\bv||_2.
\end{align}
Note that Lipschitz continuity of $\dot{\bu}(t)$ stated in \autoref{thm:existence} implies $\lim_{t\rightarrow 0^{\pm}}\partial^2_{tt} \bu(\bx, t)=\partial^2_{tt}\bu(\bx,0)$.
We now demonstrate that $\bv(\bx, t)=\partial_t \bu(\bx, t)$ is the unique solution of following initial boundary value problem.

\begin{theorem}\label{thm:existencev} 
\textbf{Initial value problem for $\bv(\bx,t)$}\\
Suppose the initial data and righthand side $\bb(t)$ satisfy the hypothesis of \autoref{thm:existence} and we suppose further that  $\dot{\bb}(t)$ exists and is continuous in time for $t\in I_0$ and $\sup_{t\in I_0} ||\dot{\bb}(t)||_2 < \infty$. Then $\bv(\bx, t)$ is the unique solution to the initial value problem $\bv(\bx, 0)=\bv_0(\bx)$,  $\partial_t \bv(\bx,0)=\partial^2_{tt}\bu(\bx, 0)$,
\begin{align}\label{eq:equationofmotionvv}
\rho \partial^2_{tt} \bv(\bx,t) &= Q(\bv(t); \bu(t))(\bx) + \dot{\bb}(\bx,t), \hbox{  $t\in I_0$},\hbox{    $\bx\in D$,}
\end{align}
$\bv \in C^2(I_0; W)$ and 
\begin{align}\label{eq:per equation vsubtt estt}
|| \partial^2_{tt} \bv(\bx,t)||_2 &\leq ||Q(\bv(t); \bu(t))(\bx)||_2 + ||\dot{\bb}(\bx,t)||_2.
\end{align}
\end{theorem}
\autoref{thm:higherregularity} now follows immediately from \autoref{thm:existencev} noting that $\partial_t \bu(\bx, t)=\bv(\bx,t)$
together with \autoref{eq:upperboundQtime} and \autoref{eq:per equation vsubtt estt}. The proof of \autoref{thm:existencev}  follows from the Lipschitz continuity \autoref{eq:lipschitzpropertyQ} and the Banach fixed point theorem as in \cite{MA-Brezis}. 

\section{Conclusions}
\label{s:conclusions}
In this article, we have provided a-priori error estimates for finite element approximations to nonlocal  state based peridynamic fracture models. We have shown that the convergence rate applies even over time intervals for which the material is softening over parts of the computational domain. The results are established for two different classes of state-based peridynamic forces. 
The convergence rate of the approximation is of the form $C(\Delta t+h^2/\epsilon^2)$ where the constant $C$ depends on $\epsilon$ and the $H^2$ norm of the solution and its time derivatives. For fixed $\Delta t$ numerical simulations for Plexiglass show that the error decreases at the rate  of $h^2$  at $40\mu$-sec into the simulation. The simulations were carried  out in parallel using 20 threads on an workstation with single Intel Xeon processor and with 32 GB of RAM. We anticipate similar convergence rates for longer times on bigger parallel machines.

{We reiterate that the a-priori error estimates account for the possible appearance of nonlinearity anywhere in the computational domain. On the other hand numerical simulation and independent theoretical estimates show that the nonlinearity concentrates along ``fat'' cracks of finite length  and width equal to $\epsilon$, see \cite{CMPer-Lipton,CMPer-Lipton3}. Moreover the remainder of the computational domain is seen to behave linearly and to leading order can be modeled as a linear elastic material up to an error proportional to $\epsilon$, see [Proposition 6, \cite{Jha and LiptonL}]. Future work will use these observations to focus on adaptive implementation and a-posteriori estimates.}


\newcommand{\noopsort}[1]{}


\begin{thebibliography}{}

\bibitem[Agwai et~al., 2011]{CMPer-Agwai}
Agwai, A., Guven, I., and Madenci, E. (2011).
\newblock Predicting crack propagation with peridynamics: a comparative study.
\newblock {\em International journal of fracture}, 171(1):65--78.

\bibitem[Aksoylu and Unlu, 2014]{AksoyluUnlu}
Aksoylu, B. and Unlu, Z. (2014).
\newblock Conditioning analysis of nonlocal integral operators in fractional
  sobolev spaces.
\newblock {\em SIAM Journal on Numerical Analysis}, 52:653--677.

\bibitem[Bobaru and Hu, 2012]{BobaruHu}
Bobaru, F. and Hu, W. (2012).
\newblock The meaning, selection, and use of the peridynamic horizon and its
  relation to crack branching in brittle materials.
\newblock {\em International journal of fracture}, 176(2):215--222.

\bibitem[Bobaru et~al., 2016]{Handbook}
Bobaru, F., Foster, J.~T., Geubelle, P.~H., Geubelle, P.~H., and Silling, S.~A.
  (2016).
\newblock Handbook of peridynamic modeling.



\bibitem[Brenner and Scott, 2007]{MANa-Susanne}
Brenner, S. and Scott, R. (2007).
\newblock The mathematical theory of finite element methods, volume~15.
\newblock {\em Springer Science \& Business Media}, 3 edition.

\bibitem[Brezis, 1983]{MA-Brezis}
Brezis, H. (1983).
\newblock Analyse fonctionnelle, th{\'e}orie et application.

\bibitem[Chen and Gunzburger, 2011]{CMPer-Chen}
Chen, X. and Gunzburger, M. (2011).
\newblock Continuous and discontinuous finite element methods for a
  peridynamics model of mechanics.
\newblock {\em Computer Methods in Applied Mechanics and Engineering},
  200(9):1237--1250.


\bibitem[Du, 2018a]{Du 2018a}
Du, Q. (2018a).
\newblock An invitation to nonlocal modeling, analysis and computation.
\newblock {\em  Proc. Int. Cong. of Math 2018}, Rio de Janeiro, Vol. 3 (3523-3552).

\bibitem[Du, 2018b]{Du 2018b}
Du, Q. (2018b).
\newblock Nonlocal modeling, analysis and computation.
\newblock NSF-CBMS Monograph, SIAM Philadelphia (2018).




\bibitem[Du et~al., 2013a]{CMPer-Du}
Du, Q., Gunzburger, M., Lehoucq, R., and Zhou, K. (2013a).
\newblock Analysis of the volume-constrained peridynamic navier equation of
  linear elasticity.
\newblock {\em Journal of Elasticity}, 113(2):193--217.

\bibitem[Du et~al., 2013b]{CMPer-Du7}
Du, Q., Ju, L., Tian, L., and Zhou, K. (2013b).
\newblock A posteriori error analysis of finite element method for linear
  nonlocal diffusion and peridynamic models.
\newblock {\em Mathematics of computation}, 82(284):1889--1922.

\bibitem[Du et~al., 2013c]{CMPer-Du6}
Du, Q., Tian, L., and Zhao, X. (2013c).
\newblock A convergent adaptive finite element algorithm for nonlocal diffusion
  and peridynamic models.
\newblock {\em SIAM Journal on Numerical Analysis}, 51(2):1211--1234.


\bibitem[Emmrich et~al., 2013]{CMPer-Emmrich}
Emmrich, E., Lehoucq, R.~B., and Puhst, D. (2013).
\newblock Peridynamics: a nonlocal continuum theory.
\newblock In {\em Meshfree Methods for Partial Differential Equations VI},
  pages 45--65. Springer.


\bibitem[Foster et~al., 2011]{CMPer-Silling7}
Foster, J.~T., Silling, S.~A., and Chen, W. (2011).
\newblock An energy based failure criterion for use with peridynamic states.
\newblock {\em International Journal for Multiscale Computational Engineering},
  9(6).

\bibitem[Demengel, 2012]{MAFa-Demengel}
Demengel, F. (2012).
\newblock Functional Spaces for the Theory of Elliptic Partial
  Differential Equations.
\newblock {\em Universitext. Springer-Verlag London}, 1 edition.


\bibitem[Gerstle et~al., 2007]{CMPer-Gerstle}
Gerstle, W., Sau, N., and Silling, S. (2007).
\newblock Peridynamic modeling of concrete structures.
\newblock {\em Nuclear engineering and design}, 237(12):1250--1258.

\bibitem[Ghajari et~al., 2014]{CMPer-Ghajari}
Ghajari, M., Iannucci, L., and Curtis, P. (2014).
\newblock A peridynamic material model for the analysis of dynamic crack
  propagation in orthotropic media.
\newblock {\em Computer Methods in Applied Mechanics and Engineering},
  276:431--452.

\bibitem[Guan and Gunzburger, 2015]{CMPer-Guan}
Guan, Q. and Gunzburger, M. (2015).
\newblock Stability and accuracy of time-stepping schemes and dispersion
  relations for a nonlocal wave equation.
\newblock {\em Numerical Methods for Partial Differential Equations},
  31(2):500--516.

\bibitem[Ha and Bobaru, 2010]{HaBobaru}
Ha, Y.~D. and Bobaru, F. (2010).
\newblock Studies of dynamic crack propagation and crack branching with
  peridynamics.
\newblock {\em International Journal of Fracture}, 162(1-2):229--244.

\bibitem[Jha and Lipton 2018b]{Jha and LiptonL}
Jha, P.~K. and Lipton, R. (2018b).
\newblock Numerical convergence of nonlinear nonlocal continuum models to local elastodynamics. 
\newblock {\em International Journal for Numerical Methods in Engineering ,}114(13):1389--1410.

\bibitem[Jha and Lipton, 2018a]{CMPer-JhaLiptonFD}
Jha, P.~K. and Lipton, R. (2018a).
\newblock Numerical analysis of nonlocal fracture models in H\"older space.
\newblock {\em SIAM J. Numer. Anal.}, 56: 906--941.

\bibitem[Jha and Lipton, 2019]{Jha and Lipton 2019}
Jha, P.~K. and Lipton, R. (2019).
\newblock Small horizon limit of state based peridynamic models.
\newblock {\em In preparation}.

\bibitem[Jha and Lipton, 2017]{CMPer-JhaLipton3}
Jha, P.~K. and Lipton, R. (2017).
\newblock Finite element approximation of nonlinear nonlocal models.
\newblock {\em arXiv preprint arXiv:1710.07661}.

\bibitem[Karaa, 2012]{CMPer-Karaa}
Karaa, S. (2012).
\newblock Stability and convergence of fully discrete finite element schemes
  for the acoustic wave equation.
\newblock {\em Journal of Applied Mathematics and Computing}, 40(1-2):659--682.


\bibitem[Lipton, 2014]{CMPer-Lipton3}
Lipton, R. (2014)
\newblock Dynamic brittle fracture as a small horizon limit of
  peridynamics. 
 \newblock {\em Journal of Elasticity}, 117:21--50.


\bibitem[Lipton, 2016]{CMPer-Lipton}
Lipton, R. (2016).
\newblock Cohesive dynamics and brittle fracture.
\newblock {\em Journal of Elasticity}, 124(2):143--191.

\bibitem[Lipton et~al., 2018a]{CMPer-Lipton5}
Lipton, R., Said, E., and Jha, P. (2018a).
\newblock Free damage propagation with memory.
\newblock {\em Journal of Elasticity}, 133(2):129--153.

\bibitem[Lipton et~al., 2018b]{CMPer-Lipton4}
Lipton, R., Said, E., and Jha, P.~K. (2018b).
\newblock Dynamic brittle fracture from nonlocal double-well potentials: A
  state-based model.
\newblock {\em Handbook of Nonlocal Continuum Mechanics for Materials and
  Structures}, pages 1--27.

\bibitem[Lipton et~al., 2016]{CMPer-Lipton2}
Lipton, R., Silling, S., and Lehoucq, R. (2016).
\newblock Complex fracture nucleation and evolution with nonlocal
  elastodynamics.
\newblock {\em arXiv preprint arXiv:1602.00247}.

\bibitem[Littlewood, 2010]{CMPer-Littlewood}
Littlewood, D.~J. (2010).
\newblock Simulation of dynamic fracture using peridynamics, finite element
  modeling, and contact.
\newblock In {\em Proceedings of the ASME 2010 International Mechanical
  Engineering Congress and Exposition (IMECE)}.

\bibitem[Macek and Silling, 2007]{CMPer-Richard}
Macek, R.~W. and Silling, S.~A. (2007).
\newblock Peridynamics via finite element analysis.
\newblock {\em Finite Elements in Analysis and Design}, 43(15):1169--1178.

\bibitem[Mengesha and Du, 2015]{CMPer-Mengesha2}
Mengesha, T. and Du, Q. (2015).
\newblock On the variational limit of a class of nonlocal functionals related
  to peridynamics.
\newblock {\em Nonlinearity}, 28(11):3999.

\bibitem[Ren et~al., 2017]{CMPer-Ren}
Ren, B., Wu, C., and Askari, E. (2017).
\newblock A 3d discontinuous galerkin finite element method with the bond-based
  peridynamics model for dynamic brittle failure analysis.
\newblock {\em International Journal of Impact Engineering}, 99:14--25.

\bibitem[Silling et~al., 2010]{CMPer-Silling5}
Silling, S., Weckner, O., Askari, E., and Bobaru, F. (2010).
\newblock Crack nucleation in a peridynamic solid.
\newblock {\em International Journal of Fracture}, 162(1-2):219--227.

\bibitem[Silling, 2000]{CMPer-Silling}
Silling, S.~A. (2000).
\newblock Reformulation of elasticity theory for discontinuities and long-range
  forces.
\newblock {\em Journal of the Mechanics and Physics of Solids}, 48(1):175--209.

\bibitem[Silling and Bobaru, 2005]{SillBob}
Silling, S.~A. and Bobaru, F. (2005).
\newblock Peridynamic modeling of membranes and fibers.
\newblock {\em International Journal of Non-Linear Mechanics}, 40(2):395--409.

\bibitem[Silling et~al., 2007]{States}
Silling, S.~A., Epton, M., Weckner, O., Xu, J., and Askari, E. (2007).
\newblock Peridynamic states and constitutive modeling.
\newblock {\em Journal of Elasticity}, 88(2):151--184.

\bibitem[Silling and Lehoucq, 2008]{CMPer-Silling4}
Silling, S.~A. and Lehoucq, R.~B. (2008).
\newblock Convergence of peridynamics to classical elasticity theory.
\newblock {\em Journal of Elasticity}, 93(1):13--37.

\bibitem[Tian and Du, 2014]{Tian-Du}
Tian, X. and Du, Q.
\newblock	Asymptotically compatible schemes and applications to robust discretization of nonlocal models.
\newblock{\em SIAM J. Numer. Anal.}, 52-4 (2014), pp. 1641--1665.

\bibitem[Weckner and Abeyaratne, 2005]{WeckAbe}
Weckner, O. and Abeyaratne, R. (2005).
\newblock The effect of long-range forces on the dynamics of a bar.
\newblock {\em Journal of the Mechanics and Physics of Solids}, 53(3):705--728.

\end{thebibliography}
\end{document}